\input amstex
\font\tenscr=rsfs10
\font\sevenscr=rsfs7
\font\fivescr=rsfs5
\skewchar\tenscr='177 \skewchar\sevenscr='177 \skewchar\fivescr='177
\newfam\scrfam \textfont\scrfam=\tenscr \scriptfont\scrfam=\sevenscr
\scriptscriptfont\scrfam=\fivescr
\define\scr#1{{\fam\scrfam#1}}
\input xy
\xyoption{all}
\loadbold
\documentstyle{amsppt}
\magnification 1100
\NoBlackBoxes
\pagewidth{450pt}

\define\Ph{\varphi}
\define\G{\text{\rm Gal}}\define\Dc{\bold D_{\text{\rm cris}}}\define\La{\varLambda}
\define\Dd{\bold D_{\text{\rm dR}}}

\define\Dpst{\bold D_{\text{\rm pst}}}
\define\Drig{\bold D_{\text{\rm rig}}}

\define\res{\text{\rm res}}

\define \bD{\bold D}
\define\rg{\text{\rm rg}}

\define \cl{\text{\rm cl}}

\define \Iw{\text{\rm Iw}}

\define\ep{\varepsilon}

\define\Bc{\bold B_{\text{\rm cris}}}\define\Exp{\text{\rm Exp}}

\define\M {\Cal M}

\define\F{\text{\rm Fil}}
\define\Fr{\text{\rm Fr}}
\define \g{\gamma}
\define \gn{\gamma_n}

\define\HG{\scr H(\Gamma)}

\define\Hi{H_{\text{\rm Iw}}}

\define\Ddag{\bold D^{\dagger}}
\define\Ddagrig{\bold D^{\dagger}_{\text{\rm rig}}}



\define\Gal{\text{\rm Gal}}
\define\dop{\partial}

\define \Dst{\bold D_{\text{\rm st}}}
\define \CDcris{\scr D_{\text{\rm cris}}}
\define  \CDst{\scr D_{\text{\rm st}}}
\define \CDpst{\scr D_{\text{\rm pst}}}
\define  \CDdr{\scr D_{\text{\rm dR}}}

\define \CE{\scr E}

\define \CR{\scr R}

\define \ea{e_\alpha}
\define \eb{e_\beta}
\define \Tw{\text{\rm Tw}}
\define \Log{\frak L}

\title
\nofrills
 Trivial zeros of $p$-adic $L$-functions at near central points\\
(second version)
\endtitle
\author
Denis Benois 
\endauthor
\abstract Using the  $\ell$-invariant constructed in our previous paper  we prove a
Mazur-Tate-Teitelbaum style formula for derivatives of $p$-adic $L$-functions of modular forms at
 near central points.
\endabstract
\leftheadtext{}
\address
Institut de Math\'ematiques de Bordeaux, UMR 5251, Universit\'e de 
Bordeaux, 351 cours de la Lib\'eration, F-33400
Talence, France 
\endaddress
\date September 5th 2012
\enddate
\email 
denis.benois\@math.u-bordeaux1.fr
\endemail
\toc \nofrills {\bf Table of contents}
\widestnumber \head{\S 5.} 
\widestnumber\subhead{5.2}
\head \S0. Introduction
\endhead
\subhead 0.1. Trivial zeros of modular forms
\endsubhead
\subhead 0.2. The semistable case
\endsubhead
\subhead 0.3. The general case
\endsubhead
\subhead 0.4. Trivial zeros of Dirichlet $L$-functions
\endsubhead
\subhead 0.5. The plan of the paper
\endsubhead
\subhead Acknowledgements
\endsubhead 
\head \S1. Preliminaries
\subhead 1.1. $(\Ph,\Gamma)$-modules
\endsubhead
\subhead 1.2. Cohomology of $(\Ph,\Gamma)$-modules
\endsubhead
\subhead 1.3. The large exponential map
\endsubhead
\subhead 1.4. $p$-adic distributions
\endsubhead
\endhead 
\head \S2. The $\ell$-invariant
\endhead
\subhead 2.1. The $\ell$-invariant
\endsubhead
\subhead 2.2. Relation to the large exponential map
\endsubhead
\head \S3. Trivial zeros of Dirichlet $L$-functions
\endhead
\subhead 3.1. Dirichlet $L$-functions
\endsubhead
\subhead 3.2. $p$-adic representations associated to Dirichlet characters
\endsubhead
\subhead 3.3. Trivial zeros 
\endsubhead
\head \S4. Trivial zeros of modular forms 
\endhead
\subhead 4.1.  $p$-adic $L$-functions 
\endsubhead
\subhead 4.2. Selmer groups and  $\ell$-invariants of modular forms
\endsubhead
\subhead 4.3. The main result
\endsubhead
\endtoc
\nologo
\endtopmatter
\document

\head
{\bf \S 0. Introduction}
\endhead

{\bf 0.1. Trivial zeros of modular forms.} In this paper we prove a  Mazur-Tate-Teitelbaum style formula for the values
of  derivatives of  $p$-adic $L$-functions of  modular forms at near central points. Together with the
results of Kato-Kurihara-Tsuji  and Greenberg-Stevens on the Mazur-Tate-Teitelbaum conjecture
this gives a complete proof of the trivial zero conjecture  formulated in \cite{Ben2} for elliptic modular forms.   
 Namely, let $f=\underset{n=1}\to{\overset \infty \to \sum} a_nq^n$ be a normalized newform on $\Gamma_0(N)$ of weight 
$k\geqslant 2$
and character $\ep$ and let $L(f,s)=\underset{n=1}\to{\overset \infty \to \sum}a_nn^{-s}$  be the complex $L$-function 
associated to $f$. It is well known that $L(f,s)$ converges for $\text{\rm Re}(s)> \frac{k+1}{2}$ and 
decomposes into an Euler product
$$
L(f,s)=\underset{l}\to \prod E_l(f,l^{-s})^{-1}
$$
where $l$ runs over all primes and  $E_l(f,X)=1-a_lX+\ep (l)l^{k-1}X^2.$ Moreover $L(f,s)$ has an analytic continuation on the whole complex plane
and satisfies the functional equation 
$$
(2\pi)^{-s}\,\Gamma (s)\, L(f,s)= i^kc N^{k/2-s}(2\pi)^{s-k}\,\Gamma (k-s)\,L(f^*,k-s)
$$
where $f^*=\underset{n=1}\to{\overset \infty \to \sum}\bar a_n q^n$ is the dual cusp form and $c$ is 
some constant (see for example \cite{Mi, Theorems 4.3.12 and 4.6.15}).
More generally, to any  Dirichlet character $\eta$   we can associate the  $L$-function
$$
L(f,\eta,s)=\underset{n=1}\to{\overset \infty \to \sum} \frac{\eta(n)\, a_n}{n^s}.
$$
The theory of modular symbols implies that  there exist non-zero complex numbers $\Omega_f^{+}$ and $\Omega_f^{-}$ such that
for any Dirichlet character $\eta$ one has 
$$
\widetilde L(f,\eta,j)=\frac{\Gamma (j)}{(2\pi i)^{j-1}\,\Omega_f^{\pm}}\,\, L(f,\eta,j) \in \overline{ \Bbb Q}, \qquad 1\leqslant j\leqslant k-1 \tag{1}
$$
where $\pm =(-1)^{j-1} \eta (-1).$ Fix a prime number $p>2$ such that the Euler factor $E_p(f,X)$ is not equal to $1$. 
Let $\alpha$ be a root of the polynomial $X^2-a_pX+\ep (p)p^{k-1}$ in $\overline {\Bbb Q}_p.$ 
Denote by $v_p$ the $p$-adic valuation on $\overline {\Bbb Q}_p$normalized such that $v_p(p)=1.$
Assume that $\alpha$ is  not critical i.e. that   $v_p(\alpha)<k-1.$ Let $\omega\,:\,(\Bbb Z/p\Bbb Z)^*@>>>\Bbb Q_p^*$
denote the Teichm\"uller character. Manin \cite{Mn}, Vishik  \cite{Vi} and independently Amice-Velu \cite{AV}  constructed 
analytic $p$-adic $L$-functions $L_{p,\alpha}(f,\omega^m,s)$ which interpolate algebraic parts of
special values of $L(f,s)$
\footnote{This construction was recently generalized to the critical case by Pollack-Stevens [PS]
and Bellaiche \cite{Bel}}. Namely, the interpolation property writes  
$$
L_{p,\alpha}(f, \omega^m,j)= \Cal E_{\alpha} (f,\omega^m,j)\,\widetilde L(f, \omega^{j-m},j), \qquad 1\leqslant j
\leqslant k-1
$$
where  $\Cal E_{\alpha} (f,\omega^m,j)$ is an explicit Euler like factor.
One says that $L_{p,\alpha}(f, \omega^m,s)$ has a trivial zero at $s=j$ if 
$\Cal E_{\alpha} (f,\omega^m,j)=0$ and $\widetilde L(f, \omega^{j-m},j)\ne 0.$ This phenomenon was first studied by Mazur, Tate and Teitelbaum 
in \cite{MTT} where the following cases were distinguished:
\newline
\newline
$\bullet$ {\it The semistable case: $p\parallel N,$ $k$ is even and $\alpha=a_p=p^{k/2-1}.$}
The $p$-adic $L$-function $L_{p,\alpha}(f, \omega^{k/2},s)$ has a trivial zero at the central point $s=k/2.$
\newline
\newline
 $\bullet${\it The crystalline case: $p\nmid N$, $k$ is odd and either $\alpha=p^{\frac{k-1}{2}}$ or 
 $\alpha= \ep (p)\,p^{\frac{k-1}{2}}$.} The $p$-adic $L$-function $L_{p,\alpha}(f, \omega^{\frac{k+1}{2}},s)$ (respectively 
$L_{p,\alpha}(f, \omega^{\frac{k-1}{2}},s)$)
has a trivial zero at the near central point
$s=\frac{k+1}{2}$ (respectively $s=\frac{k-1}{2}$).
\newline
\newline
  $\bullet${\it The potentially crystalline case: $p\mid N$, $k$ is odd  and $\alpha=a_p=p^{\frac{k-1}{2}}.$}
The  $p$-adic $L$-function $L_{p,\alpha}(f,\omega^{\frac{k+1}{2}},s)$ has a trivial zero at the near central point $s=\frac{k+1}{2}.$ 
\newline
\newline
{\bf 0.2. The semistable case.} Let 
$$
\rho_{f}\,:\,\text{\rm Gal} (\overline{\Bbb Q}/\Bbb Q) @>>> \text{\rm GL}(W_{f}).
$$
be the $p$-adic Galois representation associated to $f$ by Deligne \cite{D1}. 
Assume that $k$ is even, $p\parallel N$ and $a_p=p^{k/2-1}.$ Then the restriction of $\rho_f$ 
on the decomposition group at $p$ is semistable and non-crystalline in the sense of Fontaine \cite{Fo3}.
The associated filtered $(\Ph,N)$-module $\Dst (W_f)$ 
has a basis $\ea, \eb$ such that $\ea=N\eb,$ 
$\Ph (\ea)= a_p \ea$ and $\Ph (\eb)=p\,a_p\eb.$ The jumps of the canonical decreasing filtration of $\Dst (W_f)$ are $0$ and $k-1$
and the $\scr L$-invariant of Fontaine-Mazur is defined to be the unique element
$\scr L_{\text{FM}}(f)\in {\overline{\Bbb Q}}_p$ such that $\F^{k-1}\Dst (V_f)$ is generated by $\eb+\scr L_{\text{FM}}(f)\ea.$
In \cite{MTT} Mazur, Tate and Teitelbaum conjectured that
$$
L_{p,\alpha}(f,\omega^{k/2},k/2 )=\scr L_{\text{FM}}(f)\,\widetilde  L(f,k/2). \tag{2}
$$ 
We remark that $L(f,k/2)$ can vanish. The conjecture (2) was  proved in \cite{GS} in the weight two case
and in \cite{St} in general using the theory of $p$-adic families of modular forms. 
Another proof, based on the theory of Euler systems was found
by Kato, Kurihara and Tsuji (unpublished but see \cite{Ka2}, \cite{PR5}, \cite{Cz3}).
Note that in \cite{St} Stevens uses another definition of the $\scr L$-invariant 
proposed by Coleman \cite{Co}. 
We refer to \cite{CI} and  to the survey article \cite{Cz4} for further information and references. 
\newline
\newline
{\bf 0.3. The general case.} Our main aim in this paper is to prove an analogue of the formula (2) in the crystalline and potentially crystalline 
cases. In fact, we will treat all three cases  simultaneously.  
Let  $f$ be a newform  of  weight $k$. Fix an odd  prime   $p$ and assume that 
the $p$-adic $L$-function $L_{p,\alpha}(f, \omega^{k_0},s)$ has a trivial zero at $s=k_0.$
Assume further that $k_0\geqslant k/2.$ Note that the last assumption holds automatically in the semistable and potentially crystalline
cases and in the crystalline case it is not restrictive because we can use  the functional equation (see Remark 3 below). 
Let $\Dc (W_f)$ denote the crystalline  module associated to $W_f$. 
Our assumptions   imply that $\Dc \bigl (W_f(k_0)\bigr )^{\Ph=p^{-1}}$ is a one-dimensional vector space which we denote
by $D_{\alpha}$. The main construction of \cite{Ben2} associates to $D_{\alpha}$ an element
$\ell \left (W_f(k_0),D_{\alpha}\right )\in \overline {\Bbb Q}_p$ which can be viewed as a direct generalization of 
Greenberg's $\ell$-invariant \cite{Gre1} to the non-ordinary  case
\footnote{Strictly speaking, in \cite{Ben2} we define the $\ell$-invariant for $p$-adic representations which are semistable
at $p$, but this construction can be generalized easily to cover the potentially crystalline case (see Section 2.1 below). }. 
To simplify notation we 
set  $\ell_\alpha(f)=\ell \left (W_f(k_0),D_{\alpha}\right ).$ 
The main result of this paper states as follows.

\proclaim{\bf Theorem} Let $f$ be a newform on $\Gamma_0(N)$ of character $\ep$ and  weight $k$ and
let  $p$ be an odd prime.
Assume that  the $p$-adic $L$-function 
$L_{p,\alpha}\left (f,\omega^{k_0},s \right )$ has a trivial zero at $s=k_0\geqslant k/2$.
Then 
$$
L'_{p,\alpha}\left (f, \omega^{k_0},k_0\right )=\ell_\alpha (f)\,\left (1-\frac{\ep (p)}{p}
\right ) \widetilde L \left (f,k_0\right ).
$$ 
\endproclaim 
\flushpar 
{\bf Remarks.} 1) In the  semistable case $\ell _\alpha(f)=\scr L_{\text{FM}}(f)$ 
(see \cite{Ben2, Proposition 2.3.7}), $\ep (p)=0$ and we recover the Mazur-Tate-Teitelbaum conjecture.
Our proof in this case can be seen as an interpretation of Kato-Kurihara-Tsuji's approach in terms
of $(\Ph,\Gamma)$-modules. In the  crystalline case some version of our formula was proved by Orton \cite{Or} (unpublished). 
She does not use $(\Ph,\Gamma)$-modules and works with  an {\it ad hoc} definition of the $\ell$-invariant
in terms of Bloch-Kato exponential map in the spirit of \cite{PR4}.

2) In the crystalline and potentially crystalline cases $\widetilde L \left (f,k_0\right )$
does not vanish  by the theorem of Jacquet-Shalika [JS].

3) In the crystalline case, trivial zeros at the symmetric point $s=\frac{k-1}{2}$ can be easily
studied using the functional equation for $p$-adic $L$-functions \cite{MTT, \S17}.
If $\alpha=p^{(k-1)/2}$ and therefore $L_{p,\alpha}(f,\omega^{\frac{k+1}{2}},s)$ has a trivial zero
at $s=\frac{k+1}{2}$, then $\alpha^*= \ep^{-1} (p)\, \alpha$ is a root of the Hecke polynomial
associated to the dual form 
$f^*=\underset{n=1}\to{\overset \infty \to \sum}\bar a_nq^n$
and $L_{p,\alpha^*}(f^*,\omega^{\frac{k-1}{2}},s)$ has
a trivial zero at  $s=\frac{k-1}{2}.$  
Using the compatibility of the trivial zero conjecture with the functional equation \cite{Ben2, Section 2.3.5}
(or just repeating  the proof of the main theorem with obvious modifications)
we obtain a  trivial zero formula for $L'_{p,\alpha^*}\left (f^*,\omega^{\frac{k-1}{2}},\frac{k-1}{2}\right )$

4) The $\scr L$-invariant of Fontaine-Mazur which appears in  semistable case (2) is local i.e. it depends only on the
restriction of the $p$-adic representation $\rho_f$ on the decomposition group at $p$. 
However, in the crystalline and potentially crystalline cases our $\ell$-invariant  is global and 
contains information about the localisation map $H^1(\Bbb Q,W_f(\frac{k+1}{2}))@>>>H^1(\Bbb Q_p,W_f(\frac{k+1}{2})).$ 
We remark that in the semistable case the $p$-adic $L$-function has a zero at  the central point  and in the crystalline and potentially crystalline cases it has a zero at a near central point.

5) Let $\eta$ be a Dirichlet character of conductor $M$ with $(p,M)=1.$ The study of trivial zeros 
of $L_{p,\alpha}\left (f,\eta\omega^{\frac{k+1}{2}},s\right )$ reduces to our situation by considering
the newform $f\otimes \eta$ associated to $f_\eta=\underset{n=1}\to{\overset \infty\to\sum}\eta (n)\,a_nq^n$
(see Section 4.1.3).
\newline
\newline
Our theorem follows from a formula for  the derivative of Perrin-Riou exponential map \cite{PR2}
in terms of the $\ell$-invariant which we prove in Propositions 2.2.2 and 2.2.4 below  applied
to the Euler system constructed by Kato \cite{Ka}.
\newline
\newline
{\bf 0.4. Trivial zeros of Dirichlet $L$-functions.} Let $\eta$ be a primitive Dirichet character modulo $N$ and let $p\nmid N$ be a fixed prime. The $p$-adic $L$-function of Kubota-Leopoldt
$L_p(\eta \omega,s)$ satisfies the interpolation property
$$
L_p(\eta \omega,1-j)\,=\,(1-(\eta \omega^{1-j})(p)\,p^{j-1})\,L(\eta \omega^{1-j},1-j),\qquad j\geqslant 1.
$$ 
Assume that $\eta$ is odd and  $\eta (p)=1.$ Then $L(\eta ,0)\neq 0$ but the Euler like factor $1-(\eta \omega^{1-j})\,(p)p^{j-1}$ vanishes at $j=1$ and $L_p(\eta \omega,s)$ has a trivial zero at $s=0.$ Fix a finite extension $L/\Bbb Q_p$
containing the values of $\eta.$ Let $\chi$ denote the cyclotomic character and let $\text{\rm ord}_p\,:\,
\Gal (\Bbb Q_p^{\text{\rm ur}}/\Bbb Q_p)@>>>L$ be the character defined by $\text{\rm ord}_p(\text{\rm Fr}_p)=-1$
where $\text{\rm Fr}_p$ is the geometric Frobenius.  
Then $H^1(\Bbb Q_p,L)=\text{\rm Hom} (\Gal (\Bbb Q_p^{\text{\rm ab}}/\Bbb Q_p),L)$ is the 
two-dimensional $L$-vector space generated by $\log \chi$  and $\text{\rm ord}_p.$ Since $p\nmid N$ and $\eta (p)=1$
the restriction of $L(\eta)$ on the decomposition group at $p$ is a trivial representation. The localization map
$$
\kappa_\eta\,:\,H^1(\Bbb Q,L(\eta))@>>>H^1(\Bbb Q_p,L)
$$
is injective and identifies $H^1(\Bbb Q,L(\eta))$ with a one-dimensional subspace of $H^1(\Bbb Q_p,L).$  
It can be shown that 
$\text{\rm Im} (\kappa_\eta) $ is generated by an element of the form 
$$
\log \chi +\scr L(\eta)\,\text{\rm ord}_p \tag{3}
$$ 
there $\scr L(\eta)\in L$ is necessarily unique.
Applying Proposition 2.2.4 to the Euler system of cyclotomic units we obtain a new proof of  the trivial
zero conjecture for Dirichlet $L$-functions
$$
L'_p(\eta \omega,0)=-\scr L(\eta)\,L(\eta,0). \tag{4}
$$
This formula was first proved in \cite{Gro} as the combination of the  result of Ferrero-Greenberg \cite{FG} 
giving an explicit formula for  $L'_p(\eta\omega,0)$ in terms of the $p$-adic $\Gamma$-function and 
the Gross-Koblitz formula \cite{GK}.  We also remark that  Dasgupta, Darmon and Pollack \cite{DDP} recently 
generalized (4) to totally real number fields $F$ assuming Leopoldt's conjecture and some additional condition on the vanishing
of $p$-adic $L$-functions.
\newline
\newline
{\bf 0.5. The plan of the paper.} The main contents of this article is as follows. In \S1 we review the necessary
preliminaries. In particular, Sections 1.1-1.2 are devoted to the theory of $(\Ph,\Gamma)$-modules which plays 
a key role in our definition of the $\ell$-invariant. In Section 1.3 we review the construction and main properties
of Perrin-Riou's large exponential map.  
In \S2 we review the construction of the $\ell$-invariant $\ell (V,D)$ from \cite{Ben2}  and prove an explicit formula
for the derivative of the large logarithmic map in terms of $\ell (V,D)$  and the dual exponential map.
In \S3 we apply this formula to Dirichlet $L$-functions and give a new proof of (4). 
Trival zeros of modular forms are studied in \S4. In Section 4.1 we review basic results about $p$-adic $L$-functions of modular
forms and the representation $W_f.$ In Section 4.2 we  specialize the general definition of the $\ell$-invariant to the  case of modular forms.  Finally in Section 4.3  we  prove  the  main theorem.
\newline
\newline
{\bf Acknowledgements.} I am  grateful to Pierre Parent for a number of very valuable discussions.
A part of this work was done during my stay at the Max-Planck-Institut f\"ur Mathematik from March to May 2012. 
I would like to thank the institut for this invitation  and excellent working conditions.

\head
{\bf \S1. Preliminaries}
\endhead

{\bf 1.1. $(\Ph,\Gamma)$-modules.} 

{\bf 1.1.1. Definition of $(\Ph,\Gamma)$-modules} (see \cite{Fo1}, \cite{CC1}, \cite{Cz5}).  Let $\overline{\Bbb Q}_p$ be a fixed algebraic closure of $\Bbb Q_p$. 
We denote by $C$ the $p$-adic completion of $\overline{\Bbb Q}_p$ and $v_p\,\,:\,\,C@>>>\Bbb R\cup\{\infty\}$ the $p$-adic valuation normalized so that $v_p(p)=1$ and 
set $\vert x\vert_p=\left (\frac{1}{p}\right )^{v_p(x)}.$ Write $B(r,1)$ for the $p$-adic annulus 
$B(r,1)=\{ x\in C \,\mid \, r\leqslant \vert x\vert_p <1 \}.$
Fix a system of primitive roots of unity   $\ep=(\zeta_{p^n})_{n\geqslant 0}$,
 such that $\zeta_{p^n}^p=\zeta_{p^{n-1}}$ for all $n$. If $K$ is a finite extension of $\Bbb Q_p$
we set $G_K=\Gal (\overline{\Bbb Q}_p/K),$ $K_n= K(\zeta_{p^n})$ and $K_{\infty}= \bigcup_{n=0}^{\infty}K_n.$
Put  $\Gamma_K =\G(K_{\infty}/K)$ and denote by $\chi \,:\,\Gamma_K @>>>\Bbb Z_p^*$ the cyclotomic character.
We write $O_K$ for the ring of integers of $K$, $K_0$ for the maximal unramified subextension of $K$ and $\sigma$ for the absolute
Frobenius of $K_0/\Bbb Q_p$.

For any $r>0$ let $\bold B_K^{\dag,r}$ denote the ring of overconvergent elements of Fontaine's ring $\bold B_K$
(see \cite{CC1}, \cite{Ber1}). Note that $\bold B_K^{\dag,r_1}\subset \bold B_K^{\dag,r_2}$ if $r_1\leqslant r_2.$
The ring $\bold B_K$ is equipped with a continuous action of $\Gamma_K$ and a Frobenius operator $\Ph$
which commute to each other. We remark that  $\bold B_K^{\dag,r}$ are stable under the action of  $\Gamma_K$ 
and that $\Ph (\bold B_K^{\dag,r})\subset \bold B_K^{\dag,pr}$ for all $r$. The following description
of $\bold B_K^{\dag,r}$ is sufficent for the goals of this paper. Let $F$ denote the maximal unramified
subextension of $K_\infty/K_0$ and let $e=[K_\infty:K_0(\zeta_{p^\infty})].$
For any $0\leqslant s<1$ define
$$
\align
&\CR^{(s)}(K)\,=\,\left \{ f(X_K)=\sum_{k\in \Bbb Z}
a_kX_K^k\,\mid \, \text{\rm $a_k\in F$ and $f$ is holomorphic  on $B(s,1)$} \right \},\\
&\CE^{(s)}(K)=\,\left \{ f(X_K)=\sum_{k\in \Bbb Z}
a_kX_K^k\,\mid \, \text{\rm $a_k\in F$ and $f$ is holomorphic and bounded on $B(s,1)$} \right \}.
\endalign
$$ 
Then there exists  $r(K)>0$ such that for all $r>r(K)$ the ring $\bold B_K^{\dag,r}$ is isomorphic
to $\CE^{(p^{-1/er})}(K).$ Here the group $\Gamma_K$ acts trivially on the coefficients of power series
and $\Ph$ acts on $\CE^{(p^{-1/er})}(K)$ $\sigma$-semilineary.
In general, the action of $\Gamma_K$ and $\Ph$ on $X_K$ is very
complicated. The situation is more simple when $K=K_0$ i.e. $K$ is absolutely unramified. In this case 
$F=K_0$ and we will write $X$ instead $X_K.$ One has
$$
\align
&\tau f(X)=f(\tau (X)), \qquad \text{\rm where}\,\,\, \tau (X)=(1+X)^{\chi (\tau)}-1,\quad \tau \in \Gamma_K, \\
&\Ph f(X)= f^{\sigma}(\Ph(X)), \qquad \text{\rm where}\,\,\, \Ph(X)=(1+X)^{p}-1.
\endalign
$$

We come back to the general case. The union $\bold B_K^{\dag}=\underset{r>0}\to\cup \bold B_K^{\dag,r}$ is a 
field which is stable under the actions of $\Gamma_K$ and $\Ph$ and is isomorphic to 
$\CE^\dagger (K)=\underset{0\leqslant s<1}\to \cup \CE^{(s)}(K)$. The operator $\Ph$
has a left inverse given by
$$
\psi (f)=\frac{1}{p}\,\Ph^{-1} \left (\text{\rm Tr}_{\CE^\dagger (K)/\Ph(\CE^\dagger (K))}(f)\right ).
$$
If $K=K_0$ we can also write
$$
\psi (f(X))=\frac{1}{p}\,\Ph^{-1} \left (\underset{\zeta^p=1}\to\sum f(\zeta (1+X)-1)\right ).
$$

The field 
$\CE^\dagger (K)$ is  endowed with the valuation 
$$
w\left (\sum_{k\in \Bbb Z}
a_kX^k \right )=\min \{v_p(a_k) \mid k\in \Bbb Z\}
$$ 
and we denote by $\Cal O_{{\CE}^\dag(K)}$ its ring of integers.

Set $\CR (K)=\underset{0\leqslant s<1}
\to \cup \CR^{(s)}(K).$ The actions of $\Gamma_K$, $\Ph$ and $\psi$ can be extended
to  $\CR (K)$ by continuity.  If $K\subset K'$ then the natural inclusions
$\bold B_K^{\dag,r}\subset  \bold B_{K'}^{\dag,r}$ induce  embeddings
$\CE^\dag (K)\subset \CE^\dag (K')$ and $\CR (K)\subset \CR(K').$  
Let $t=\log (1+X)=\underset{k=1}\to{\overset{\infty}\to\sum} (-1)^{k+1}X^k/k\in \CR(\Bbb Q_p).$ 
Note that $\Ph (t)=p\,t$ and $\tau (t)=\chi (\tau)t$ for all $\tau\in \Gamma_K.$

In this paper we deal with $p$-adic representations with coefficients in a finite
extension $L$ of $\Bbb Q_p.$ By this reason it is convenient to set $\CE_L^\dagger (K)=\CE^\dagger (K)\otimes_{\Bbb Q_p}L,$
$\CR_L^\dagger (K)=\CR^\dagger (K)\otimes_{\Bbb Q_p}L$ and 
$\Cal O_{{\CE}_L^\dag(K)}=\Cal O_{{\CE}^\dag(K)}\otimes_{\Bbb Z_p}O_L.$

\proclaim{Definition} i) A $(\Ph,\Gamma_K)$-module over $\CE^\dagger_L(K)$ (resp. $\CR_L(K)$) is a free
$\CE^\dagger_L(K)$-module (resp. $\CR_L(K)$-module) $\bD$ of finite rank $d$ equipped with
semilinear actions of $\Gamma_K$ and $\Ph$ which commute to each other and such that the induced linear map
$\CE^\dagger_L(K)\otimes_{\Ph}\bD @>>>\bD$ (resp. $\CR_L(K)\otimes_{\Ph}\bD @>>>\bD$) is an isomorphism. 

ii) A $(\Ph,\Gamma_K)$-module $\bD$ over $\CE^\dagger_L(K)$ is said to be etale if there exists a basis
of $\bD$ such that the matrix of $\Ph$ in this basis is in $\text{\rm GL}_d(\Cal O_{{\CE_L}^\dag(K)}).$  
\endproclaim
\flushpar
If $\bD$ is a $(\Ph,\Gamma_K)$-module over $A=\CE_L^\dag(K)$ or $ \CR_L(K)$ we write $\bD^*$ for the dual 
module $\text{\rm Hom}_{A}(\bD,A)$ and $\bD(\chi)$ for the module obtained from $\bD$ by twisting 
the action of $\Gamma_K$ by the cyclotomic character. 
\newline
\newline
Let $\bold{Rep}_L (G_{K})$ be the category of $p$-adic representations  of $G_{K}$
with coefficients in $L$  i.e.
the category of finite dimensional $L$-vector spaces equipped with a continuous linear action of
$G_{K}.$ 

\proclaim\nofrills{Theorem 1.1.2}\,\,\rm{(\cite{Fo1}, \cite{CC1})}. There exists a natural functor
$
V@>>>\Ddag (V)
$
which induces an equivalence between $\bold{Rep}_L (G_{K})$ and the category of 
etale $(\Ph,\Gamma_K)$-modules over $\CE^\dag_L(K).$
\endproclaim
\flushpar
\,
\flushpar
From Kedlaya's theory it follows \cite{Cz5, Proposition 1.4 and Corollary 1.5} that the functor 
$\bD @>>>\CR_L(K)\otimes_{\CE^\dag_L(K)} \bD$ establishes an equivalence between the category
of \'etale $(\Ph,\Gamma_K)$-modules over $\CE^\dag_L(K)$ and the category of $(\Ph,\Gamma_K)$-modules over
$\CR_L(K)$ of slope $0$ in the sense of \cite{Ke}. Together  with Theorem 1.1.2 this implies that the functor
$V@>>>\Ddagrig (V) $ defined by
$\Ddagrig (V)=\CR_L(K)\otimes_{\CE^\dag_L(K)}\Ddag (V)$ induces an equivalence between the category of $p$-adic
representations and the category of  $(\Ph,\Gamma_K)$-modules over
$\CR_L(K)$ of slope $0$.
\newline
\newline
{\bf 1.1.3. Crystalline and semistable  $(\Ph,\Gamma_K)$-modules} (see \cite{Fo3}, \cite{Ber3}, \cite{Ber4}). Recall that a filtered $(\Ph,N)$-module over $K$
with coefficients in $L$ is a free $K_0\otimes_{\Bbb Q_p} L$-vector space $M$ equipped with the following structures:

$\bullet$ a $\sigma$-semilinear isomorphism $\Ph\,:\,M@>>>M$ ($\sigma$ acts trivially on $L$);

$\bullet$ a $K_0\otimes_{\Bbb Q_p} L$-linear nilpotent operator $N$ such that $N\,\Ph=p\,\Ph\, N;$

$\bullet$ an exhaustive decreasing filtration $(\F^i M_{K})_{i\in \Bbb Z}$ on $M_{K}=K\otimes_{K_0}M$ by 
$(K\otimes_{\Bbb Q_p} L)$-submodules.
\newline
If $K'/K$ is a finite Galois  extension with Galois group $G_{K'/K}$, then a
filtered $(\Ph,N,G_{K'/K})$-module is a filtered $(\Ph,N)$-module $M$ over $K'$ equipped
with a semilinear action of $G_{K'/K}$ such that the filtration  $\F^i M_{K'}$
is $G_{K'/K}$-stable. We say that $M$ is a filtered $(\Ph,N,G_K)$-module
if it is a filtered $(\Ph,N,G_{K'/K})$-module for some $K'/K.$
It is well known (see for example \cite{Fo3}) that filtered $(\Ph,N,G_K)$-modules form a tensor category
 $\bold M\bold F^{\Ph,N,G_K}_{K}$ which  is additive,
has kernels and cokernels but is not abelian. The unit object $\bold 1$ of $\bold M\bold F^{\Ph,N,G_K}_{K}$ 
is  the module $K_0\otimes_{\Bbb Q_p}L$
with the natural action of $\Ph$ and  the
filtration given by
$$
\F^i\bold 1_K=\cases K\otimes_{\Bbb Q_p}L,&\text{\rm if $i\leqslant 0,$}\\
0, &\text{\rm if $i>0.$}
\endcases
$$
A filtered $(\Ph,N)$-module can be viewed as a filtered $(\Ph,N,G_K)$-module
with the trivial action of $G_K$ and we denote by  $\bold M\bold F^{\Ph,N}_{K}$ 
the resulting subcategory. A filtered Dieudonn\'e module is an object
$M$ of $\bold M\bold F^{\Ph,N}_{K}$ such that $N=0$ on $M.$ Filtered Dieudonn\'e
modules form a full subcategory $\bold M\bold F^{\Ph}_{K}$ of $\bold M\bold F^{\Ph,N}_{K}$.

If $M$ is a filtered $(\Ph,N,G_K)$ - module of rank $1$ and $m$ is a basis vector
of $M$, then $\Ph (m)=\alpha \, m$ for some $\alpha \in L$. Set
$t_N(M)= v_p(\alpha)$ and denote by $t_H(M)$ the unique jump in the  filtration   
of $M.$ If $M$ has  rank $d\geqslant 1$, set 
$t_N(M)= t_N(\overset d\to \wedge M)$ and  $t_H(M)= t_H(\overset d\to \wedge M).$   
A filtered $(\Ph,N,G_K)$-module  $M$ is said to be weakly admissible if $t_H(M)=t_N(M)$
and  $t_H(M')\leqslant t_N(M')$  for any $(\Ph,N,G_K)$-submodule $M'$ of $M$. Weakly admissible modules 
form a subcategory of $\bold M\bold F^{\Ph,N,G_K}_{K}$ which we denote by $\bold M\bold F_{K,f}^{\Ph,N,G_K}.$

Let $\CR_{L,\log}(K)=\CR_L(K)[\log X]$ where $\log X$ is  transcendental over $\CR_L(K)$
and 
$$
\tau (\log X)= \log X+\log \left (\tau (X)/X\right ),\quad \tau \in \Gamma_K,
\qquad \Ph (\log X)=p\log X+\log \left (\Ph(X)/X^p\right ).
$$
Define a monodromy operator $N\,:\,\CR_{L,\log}(K)@>>>\CR_{L,\log}(K)$ by
$N=-\left (1-\dsize\frac{1}{p}\right )^{-1}\,\dsize\frac{d}{d\log X}.$
For any $(\Ph,\Gamma_K)$-module $\bD$  over $\CR_L(K)$ we set
$$
\CDcris (\bD)=\left (\bD \otimes_{\CR_L(K)}\CR_L(K) [1/t]\right )^{\Gamma_K},\qquad
\CDst (\bD)=\left (\bD \otimes_{\CR_L(K)}\CR_{L,\log}(K) [1/t]\right )^{\Gamma_K}.
$$
Then $\CDcris (\bD)$ (resp. $\CDst (\bD)$) is a   $K\otimes_{\Bbb Q_p}L$-module of finite rank equipped with
a natural action of $\Ph$ (resp. with natural actions of $\Ph$ and $N$).
There exists a compatible system of embeddings  $\Ph^{-m}\,:\,\CR_L(K)^{(r)}@>>>(L\otimes K_\infty) [[t]]$ 
 which allows to define  exhaustive decreasing filtrations
on $\CDcris (\bD)_{K}$ and $\CDst (\bD)_K)$ (see \cite{Ber4, proof of Theorem III.2.3} ). Moreover 
$\CDcris (\bD)=\CDst (\bD)^{N=0}$ and 
$$
\rg (\CDcris (\bD))\leqslant \rg (\CDst(\bD))\leqslant  \text{rg} (\bD).
$$
We say that $\bD$ is crystalline (resp. semistable) if
$\rg (\CDcris (\bD))= \text{rg} (\bD)$ (resp. $\rg (\CDst(\bD))=  \text{rg} (\bD)$). 

If $K'/K$ is a finite extention, then $\CR_L(K)\subset \CR_L(K')$ and we set 
$\bD_{K'}=\CR_L(K')\otimes_{\CR_L(K)}\bD.$ The projective limit
$$
\CDpst (\bD)=\underset{K'/K}\to\varinjlim \scr D_{\text{\rm st}/K'} (\bD_{K'})
$$
is a $K_0^{ur}\otimes_{\Bbb Q_p} L$-module of finite rank equipped with a discrete action of $G_K$ 
and we say that $\bD$ is potentially semistable if $\rg (\CDpst (\bD))= \text{rg} (\bD).$
Denote by $\bold M^{\,\Ph,\Gamma}_{\text{pst},K}$, ,$\bold M^{\,\Ph,\Gamma}_{\text{st},K}$
and $\bold M^{\,\Ph,\Gamma}_{\text{cris},K}$ the categories of potentially semistable,
semistable and crystalline $(\Ph,\Gamma_K)$-modules respectively.

\proclaim{Proposition 1.1.4} i) The functors 
$ \CDcris \,:\,\bold M^{\,\Ph,\Gamma}_{\text{cris},K}@>>> \bold M\bold F^{\Ph}_{K},$ 
$ \CDst \,:\,\bold M^{\,\Ph,\Gamma}_{\text{st},K}@>>> \bold M\bold F^{\Ph,N}_{K}$ and
$ \CDpst \,:\,\bold M^{\,\Ph,\Gamma}_{\text{pst},K}@>>> \bold M\bold F^{\Ph,N,G_K}_{K}$
are equivalences of categories.

ii) If $V$ is a $p$-adic representation of $G_{K}$ then $\CDcris (\Ddagrig (V))$ (resp.
$\CDst (\Ddagrig (V))$, resp.  $\CDpst (\Ddagrig (V))$)
is canonically
and fonctorially isomorphic to Fontaine's module $\Dc (V)$  (resp. $\Dst (V$), resp. $\Dpst (V)$). 
\endproclaim 
\demo{Proof} The first statement is the main result of \cite{Ber4}. The second statement
follows from \cite{Ber1, Theorem  0.2}. 
\enddemo 
$\,$
\newline
\newline
{\bf 1.2. Cohomology of $(\Ph,\Gamma)$-modules.}
\newline
{\bf 1.2.1. Fontaine-Herr complexes} (see \cite{H1}, \cite{H2}, \cite{Liu}).
Let $A$ be either $\CE_L^\dag(K)$ or $\CR_L(K).$ We fix  a generator $\g_K\in \Gamma_K $. If $\bD$ is a $(\Ph,\Gamma_K)$-module over
$A$  we shall write $H^*(\bD)$ for the cohomology of the complex
$$
C_{\Ph,\gamma_K} (\bD)\,\,:\,\,0@>>>\bD @>f>> \bD\oplus \bD@>g>> \bD@>>>0
$$
where $f(x)=((\Ph-1)\,x,(\gamma_K -1)\,x)$ and $g(y,z)=(\gamma_K-1)\,y-(\Ph-1)\,z.$
A short exact sequence of $(\Ph,\Gamma_K)$-modules
$$
0@>>>\bD'@>>>\bD@>>>\bD''@>>>0
$$
gives rise to an exact cohomology sequence:
$$
0@>>>H^0(\bD')@>>>H^0(\bD)@>>>H^0(\bD'')@>>>H^1(\bD')@>>>\cdots @>>>
H^2(\bD'')@>>>0.
$$
The cohomology of $(\Ph,\Gamma_K)$-modules over $\CR_L(K)$ satisfies the following fondamental properties
(see \cite{Liu, Theorem 0.2}):

$\bullet$ {\it Euler characteristic formula.} $H^*(\bD)$ are finite dimensional $L$-vector spaces 
and the usual formula for the Euler characteristic holds
$$
\underset{i=0}\to{\overset 2 \to\sum} (-1)^i \dim_L H^i(\bD)=-[K:\Bbb Q_p]\,\text{\rm rg}_{\CR_L(K)}(\bD).
$$ 

$\bullet$ {\it Poincar\'e duality.}  For each $i=0,1,2$ there exist functorial pairings
$$
H^i(\bD) \times H^{2-i}(\bD^*(\chi)) @>\cup >>H^2(\CR_L(K)(\chi))\simeq L
$$ 
which are compatible with the connecting homomorphisms in the usual sense.

\proclaim{Proposition 1.2.2} Let $V$ be a $p$-adic representation of $G_{K}.$
Then 

i) The continuous Galois cohomology $H^*(K,V)$ is canonically (up to the choice of $\gamma_K$) and
functorially isomorphic to $H^*(\Ddag (V)).$ 

ii) The natural map $\Ddag (V)@>>>\Ddagrig (V)$ induces a quasi-isomorphism of complexes
$C_{\Ph,\gamma_K} (\Ddag (V)) @>>> C_{\Ph,\gamma_K} (\Ddagrig (V)).$
\endproclaim
\demo{Proof} See \cite{H1} and \cite{Liu, Theorem 1.1}.
\enddemo 
\flushpar
{\bf 1.2.3. Iwasawa cohomology} (see \cite{CC2}). 
 If $V$ is a $p$-adic representation of $G_{K}$ and  $T$ is an $O_L$-lattice of  $V$ 
stable under $G_{K}$ we define
$$
H^i_{\Iw}(K,T)=\varprojlim_{\text{cor}_{K_n/K_{n-1}}} H^i(K_n,T)
$$
and $H^i_{\Iw}(K,V)=H^i_{\Iw}(K,T)\otimes_{O_L}L.$
Since  $\Ddag (V)$ is etale, each  $x\in \Ddag (V)$ can be written in the form 
$x=\dsize\underset{i=1}\to{\overset d\to\sum} a_i \Ph(e_i)$ where $\{e_i\}_{i=1}^d$ is a basis of $\Ddag (V)$
and $a_i\in \CE_L^\dag(K)$. Therefore the formula 
$$
\psi \left ( \underset{i=1}\to{\overset d\to\sum} a_i \,\Ph(e_i)\right )=\underset{i=1}\to{\overset d\to\sum} \psi (a_i)\, e_i
$$  
defines an operator $\psi \,:\,\Ddag (V)@>>>\Ddag (V)$ which is a left inverse for $\Ph.$
The Iwasawa cohomology $H^*_{\Iw}(K,V)$ is canonically (up to the choice of $\gamma_K$) and functorially isomorphic to
the cohomology of the complex
$$
C^{\dag}_{\Iw,\psi}(V)\,\,:\,\, \bD^{\dag} (V)@>\psi-1>>\bD^{\dag} (V).
$$
The projection map $\text{\rm pr}_{V,n}\,:\,H^1_{\Iw}(K,V)@>>>H^1(K_n,V)$ has the following explicit description.
Set $\gamma_{K,n}=\gamma_K^{[K_n:K]}.$
 Let $x\in \bD^{\dag}(V)^{\psi=1}.$ Then $(\Ph-1)\,x\in \bD^{\dag}(V)^{\psi=0}$ and
by  \cite{CC1, Lemma 1.5.1} there exists $y\in \bD^{\dag}(V)$ such that $(\gamma_{K,n}-1)\,y=(\Ph-1)\,x.$ Then 
$\text{\rm pr}_{V,n}$ sends $x$ to $\cl (y,x).$ 
This interpretation of the Iwasawa cohomology  was found by Fontaine 
(unpublished but see \cite{CC2}).
\newline
\newline
{\bf 1.2.4. The exponential map} (see \cite{BK}, \cite{Ne}, \cite{Ben2}). 
Let $\bD$ be a $(\Ph,\Gamma_K)$-module. To any cocycle $\alpha =(a,b)\in Z^1(C_{\Ph,\gamma}(\bD))$ 
one can associate the extension
$$
0@>>>\bD@>>>\bD_{\alpha}@>>> \CR_L(K)@>>>0
$$
defined by 
$$
\bD_{\alpha}=\bD\oplus \CR_L(K)\,e,\qquad (\Ph-1)\,e=a, \quad (\gamma_K-1)\,e=b.
$$
As usual, this gives rise to a canonical isomorphism  $H^1(\bD)\simeq \text{\rm Ext}^1_{(\Ph,\Gamma_K)} (\CR_L(K),\bD).$ 
We say that  the class $\cl (\alpha)$ of $\alpha$ in $H^1(\bD)$ is crystalline if
$\rg_{L\otimes K_0} \scr D_{\text{\rm cris}}(\bD_{\alpha})=
\rg_{L\otimes K_0} \scr D_{\text{\rm cris}}(\bD)+1
$
and define 
$$
H^1_f(\bD)\,=\,\{ \cl (\alpha)\in H^1(\bD)\,\,\vert \,\, \text{\rm $\cl (\alpha)$ is crystalline}\,\}
$$
(see \cite{Ben2, Section 1.4}). Now assume that $\bD$ is potentially semistable
and define the tangent space of $\bD$ as 
$$
t_{\bD}(K)= {\scr D_{\text{\rm dR}}(\bD)}/\F^0 {\scr D_{\text{\rm dR}}(\bD)}.
$$
Consider the complex 
$$
C^{\bullet}_{\text{\rm cris}}(\bD )\,\,:\,\,\CDcris (\bD)@>f>> t_{\bD}(K) \oplus \CDcris (\bD)
$$
where the modules are placed in degrees $0$ and $1$ and $f(d)=
(d\pmod {\F^0 \scr D_{\text{\rm dR}}(\bD)},(1-\Ph)\,(d))$ (see \cite{Ne}, \cite{FP}).
From Proposition 1.1.4 it follows the existence of a canonical isomorphism
$$
H^1( C^{\bullet}_{\text{\rm cris}}(\bD))@>>>H^1_f(\bD)
$$
(see \cite{Ben2, Proposition 1.4.4} for the proof).
We define the exponential map
$$
\exp_{\,\bD,K}\,\,:\,\,t_{\bD}(K)\oplus \CDcris (\bD)@>>>H^1(\bD)
$$
as the composition of this isomorphism with the natural projection 
$t_{\bD}(K)\oplus \CDcris (\bD) @>>>H^1(C^{\bullet}_{\text{\rm cris}}(\bD))$ and the embedding
$H^1_f(\bD)\hookrightarrow H^1(\bD).$

If $V$ is a potentially semistable  representation and $\bD=\Ddagrig (V)$ then the isomorphism $H^1(\bD)\simeq H^1(K,V)$ 
identifies $H^1_f(\bD)$ with  $H^1_f(K,V)$ of Bloch-Kato \cite{Ben2, Proposition 1.4.2}.
Let 
$$t_V(K)=\Dd (V)/\F^0\Dd(V)
$$ 
denote the tangent space of $V$.  
By \cite{Ne, Proposition 1.21} the following diagram commutes and identifies our exponential map
with the exponential map $\exp_{V,K}$ of Bloch-Kato \cite{BK, \S 4}
$$
\CD
t_{\bD}(K) @>\exp_{\,\bD,K}>> H^1(\bD)\\
@VVV @VVV\\
t_V(K) @>\exp_{V,K_n}>>H^1(K,V).
\endCD
$$
Let 
$$
\left [\,\,,\,\,\right ]\,\,:\,\, {\CDdr(\bD)} \times {\CDdr(\bD^* (\chi))} @>>>L\otimes_{\Bbb Q_p}K
$$
be the canonical duality. The dual exponential map
$$
\exp^*_{\,\bD^* (\chi),K}\,\,:\,\,H^1(\bD^*(\chi)) @>>>\F^0 {\CDdr (\bD^*(\chi))} 
$$
is defined to be the unique linear map such that
$$
\exp_{\bD,K} (x) \cup y=\text{\rm Tr}_{K/\Bbb Q_p} [x,\exp^*_{\bD^*(\chi),K}(y)]
$$
for all $x\in {\CDdr (\bD)}$, $y\in {\CDdr (\bD^* (\chi))}.$
\newline
\newline
{\bf 1.2.5. $(\Ph,\Gamma)$-modules of rank $1$} (see \cite{Cz5}, \cite{Ben2}).
In this paper we deal with potentially semistable representations of $G_{\Bbb Q_p}.$
To simplify notation we set $K_n=\Bbb Q_p(\zeta_{p^n}),$  
$\CE_L^\dag=\CE_L^\dag (\Bbb Q_p),$  $\CR_L=\CR_L(\Bbb Q_p),$ $\Gamma=\Gamma_{\Bbb Q_p}$ 
and we fix a topological generator $\gamma$ of $\Gamma.$
With each continuous character $\delta\,:\,\Bbb Q^*_p@>>>L^*$ one can  associate  the $(\Ph,\Gamma)$-module
of rank one $\CR_L(\delta)=\CR_Le_{\delta}$ defined by  $\gamma (e_\delta)=\delta (\chi (\gamma))e_\delta$
and $\Ph (e_\delta)=\delta (p)e_{\delta}$.  Colmez proved that any $(\Ph,\Gamma)$-module  
of rank one over $\CR_L$ is isomorphic to one and only one of  $\CR_L(\delta)$ \cite{Cz5, Proposition 3.1}.
It is easy to see that $\CR_L(\delta)$ is crystalline if and only if there exists $m\in \Bbb Z$ such that
$\delta (u)=u^m$ for all $u\in \Bbb Z_p^*$ \cite{Ben2, Lemma 1.5.2}. In this case $\CDcris (\CR_L(\delta))$
is the one-dimensional vector space generated by $t^{-m}e_\delta$ with Hodge-Tate weight equal to $-m$
 and $\Ph$ acts on  $\CDcris (\CR_L(\delta))$ as multiplication by $p^{-m}\delta (p).$ 
The computation of the cohomology of crystalline $(\Ph,\Gamma)$-modules of rank $1$ reduces to the following 
four cases. We refer to \cite{Cz5,  Sections 2.3-2.5}  and to \cite{Ben2, Proposition 1.5.3 and Theorem 1.5.7} for
proofs and more details.   
\newline
\,

$\bullet$ {\it $\delta (u)=u^{-m}$ {\rm (}$u\in \Bbb Z_p^*${\rm )}  for some $m\geqslant 0$ but  
$\delta (x)\neq x^{-m}$.} In this case $H^i(\CR_L(\delta))=0$ for $i=0,2$, 
$H^1(\CR_L(\delta))$ is a one-dimensional $L$-vector space and  $H^1_f(\CR_L(\delta))=0.$
\newline
\,

$\bullet$ {\it $\delta (x)=x^{-m}$ for some $m\geqslant 0$.}  In this case  $H^0(\CR_L(\delta))=\CDcris(\CR_L(\delta))$ and $H^2(\CR_L(\delta))=0.$
The map
$$
\align
&i_{\delta}\,:\, \CDcris(\CR_L(\delta))\oplus \CDcris(\CR_L(\delta)) @>>>H^1(\CR_L(\delta)),\\
&i_{\delta}(x,y)=\cl (-x,\log \chi (\gamma)y)
\endalign
$$
is an isomorphism.  We let  $i_{\delta,f}$ and $i_{\delta,c}$ denote its restrictions on the first and second
direct summand respectively. Then $\text{\rm Im} (i_{\delta,f})=H^1_f(\CR_L(\delta))$ and  we have a canonical decomposition
$$
H^1(\CR_L(\delta))\simeq H^1_f(\CR_L(\delta))\oplus H^1_c(\CR_L(\delta)) \tag{5} 
$$
where $H^1_c(\CR_L(\delta))=\text{\rm Im} (i_{\delta,c}).$  Set 
$$
\align
&{\bold x}_m=i_{\delta,f}(t^{m}e_{\delta})=-\cl (t^{m},0)\,e_\delta,\\ 
&\bold y_m=i_{\delta,c}(t^{m}e_{\delta})=\log \chi (\g)\,\cl (0,t^{m})\,e_\delta.
\endalign
$$

$\bullet$ {\it $\delta (u)=u^m$ ($u\in \Bbb Z_p^*$) for some $m\geqslant 1$ but $\delta (x)\neq \vert x\vert x^m$.} Then $H^i(\CR_L(\delta))=0$ for $i=0,2$, 
$H^1(\CR_L(\delta))$ is a one-dimensional $L$-vector space and  $H^1_f(\CR_L(\delta))=H^1(\CR_L(\delta)).$
\newline
\,

$\bullet$ {\it $\delta (x)=\vert x\vert x^m$ for some $m\geqslant 1.$}  Then $H^0(\CR_L(\delta))=0$ and 
$H^2(\CR_L(\delta))$ is a one-dimensional $L$-vector space. Moreover $\chi\delta^{-1} (x)= x^{1-m}$ 
and there exists a unique isomorphism
$$
i_{\delta}\,:\, \CDcris(\CR_L(\delta))\oplus \CDcris(\CR_L(\delta)) @>>>H^1(\CR_L(\delta))
$$
such that 
$$
i_{\delta}(\alpha,\beta) \cup i_{\chi\delta^{-1}}(x,y)=[\beta,x]-[\alpha,y]
$$
where $[\,\,,\,]\,:\,\CDcris (\CR_L(\delta))\times \CDcris (\CR_L(\chi\delta^{-1})) @>>>L$ is the canonical pairing.
 Denote $i_{\delta,f}$ and $i_{\delta,c}$ the restrictions of $i_{\delta}$  on the first and second
direct summand respectively. Then $\text{\rm Im} (i_{\delta,f})=H^1_f(\CR_L(\delta))$ and again  we have a canonical decomposition
$$
H^1(\CR_L(\delta))\simeq H^1_f(\CR_L(\delta))\oplus H^1_c(\CR_L(\delta)) \tag{6}
$$
where $H^1_c(\CR_L(\delta))=\text{\rm Im} (i_{\delta,c}).$ 

More  explicitly, let
$\boldsymbol \alpha_m =-\dsize\left (1-\frac{1}{p}\right )\,\cl (\alpha_m)$ and 
$\boldsymbol \beta_m =\dsize\left (1-\frac{1}{p}\right ) \log \chi (\gamma)\,\cl (\beta_m)$
where
$$
\align
&\alpha_m=\frac{(-1)^{m-1}}{(m-1)!}\, \partial^{m-1} \left (\frac{1}{X}+\frac{1}{2},a \right )\,e_{\delta} ,\qquad
(1-\Ph)\,a=(1-\chi (\gamma) \gamma)\,\left (\frac{1}{X}+\frac{1}{2} \right ),\\
&\beta_m=\frac{(-1)^{m-1}}{(m-1)!}\,\partial^{m-1} \left (b,\frac{1}{X} \right )\,e_{\delta}, \qquad
(1-\Ph)\,\left (\frac{1}{X}\right )\,=\,(1-\chi (\gamma)\,\gamma)\,b
\endalign
$$
and $\partial =(1+X)\dsize\frac{d}{dX}.$ Then $H^1_f(\CR_L(\delta))$ and $H^1_c(\CR_L(\delta))$ are generated by
 $\boldsymbol\alpha_m$ and $\boldsymbol\beta_m$ respectively and one has
$$
\boldsymbol\alpha_m\cup \bold x_{m-1}=\boldsymbol\beta_m\cup \bold y_{m-1}=0, 
\quad \boldsymbol \alpha_m\cup \bold y_{m-1}=-1, \quad 
\boldsymbol \beta_m\cup \bold x_{m-1}=1. \tag{7}
$$

\proclaim{Proposition 1.2.6} Let $\delta (x)=|x|x^m$ where $m\geqslant 1.$ Then
$d_m=t^{-m}e_\delta$ is a basis of $\CDcris (\CR_L(\delta))$ and 
the exponential map sends $(d_m,0)$ to $\boldsymbol\alpha_m$.
\endproclaim
\demo{Proof} See \cite{Ben2, Theorem 1.5.7}.
\enddemo
\flushpar
\, 
\flushpar
{\bf 1.3. The large exponential map.}
\newline
{\bf 1.3.1. The large exponential map} (see \cite{PR2}, \cite{Cz1}, \cite{Ben1}, \cite{Ber2}). In this section we review the construction and basic properties of 
Perrin-Riou's large exponential map \cite{PR2}.  We work with $p$-adic representations
of $G_{\Bbb Q_p}$ and keep notations of Section 1.2.5.
Let  $p$ be an odd prime number.  We let denote $\Lambda=O_L[[\Gamma]]$ the Iwasawa
algebra of $\Gamma$ over $O_L$ and set $\CR_L^+=\CR_L\cap L[[X]].$ We remark that $\CR_L^+$ is the ring
of power series with coefficients in $L$ which converge on the open unit disk. 
Fix a topological generator $\gamma$ of $\Gamma$ and define a compartible system of generators of $\Gamma_n$ setting
$\gamma_1=\gamma^{p-1}$ and $\gamma_{n+1}=\gamma_n^p$ for $n\geqslant 1.$ Let $\Delta=\text{\rm Gal}(K_1/\Bbb Q_p).$
Define 
$$
\scr H\,=\,\{f(\g_1-1)\,|\,f \in \CR_L^+ \},\qquad
\scr H(\Gamma)\,=\,\Bbb Z_p[\Delta]\otimes_{\Bbb Z_p}\scr H (\Gamma).
$$
Thus $\scr H(\Gamma)=\underset{i=0}\to{\overset{p-2}\to \oplus}\scr H \delta_i$ where
$\delta_i=\dsize\frac{1}{\vert\Delta\vert}\dsize \underset{g\in \Delta}\to \sum \omega^{-i}(g) g.$ 
We equip $\scr H(\Gamma)$ with twist operators $\Tw_m\,:\,\scr H(\Gamma) @>>>\scr H(\Gamma)$ defined by $\Tw_m(f (\gamma_1-1)\, \delta_i)=
f(\chi (\gamma_1)^m\gamma_1-1)\,\delta_{i-m}.$
The ring $\scr H (\Gamma)$ acts on $\CR_L^+$ and $(\CR_L^+)^{\psi=0}$ is the free $\scr H (\Gamma)$-module generated by $(1+X)$
\cite{PR2, Proposition 1.2.7}.
\newline
\newline
Let $V$ be a potentially semistable representation of $G_{\Bbb Q_p}.$ 
Set $\Cal D (V)=(\CR_L^+)^{\psi=0}\otimes_L \Dc (V)$ and define a map
$$
\Xi^\ep_{V,n}\,:\,\Cal D(V)@>>>H^1(K_n,C^{\bullet}_{\text{cris}}(\Ddagrig (V)))=\text{\rm coker} \left (\Dc (V)@>f>>t_V(K_n)\oplus \Dc (V)\right )
$$ 
by
$$
\Xi_{V,n}^\ep (\alpha)\,=\,\cases p^{-n} \biggl (\sum_{k=1}^n
(\sigma \otimes \Ph)^{-k} \alpha (\zeta_{p^k}-1),\,-\alpha (0) \biggr )
&\text{\rm if
$n\geqslant 1$,} \\
-\bigl (0,(1-p^{-1}\Ph^{-1})\,\alpha (0) \bigr )
&\text{\rm if $n=0.$}
\endcases
$$
In particular, if $\Dc (V)^{\Ph=1}=0$ the operator $1-\Ph$ is invertible on $\Dc (V)$
and 
$$\Xi^{\ep}_{V,0}(\alpha)=\dsize \left ({\frac{1-p^{-1}\Ph^{-1}}{1-\Ph}}\,\alpha (0),0 \right ).
$$
For any $m\in \Bbb Z$ let $\Tw_{V,m}^{\ep} \,:\,H^1_{\Iw}(\Bbb Q_p, V)@>>>H^1_{\Iw}(\Bbb Q_p,V(m))$ denote 
the twist map $\Tw_{V,m}^\ep (x)=x\otimes \ep^{\otimes m}.$

\proclaim{Theorem 1.3.2} Let $V$ be a potentially semistable  representation
of $G_{\Bbb Q_p}$ such that $H^0 (K_\infty,V)=0$. Then  for any integers 
$h$ and $m$ such that $\F^{-h}\Dd (V)=\Dd (V)$ and $m+h\geqslant 1$   there exists 
a unique $\scr H (\Gamma)$-homomorphism
$$
\Exp_{V (m),h}^\ep \,\,:\,\,\Cal D(V(m))  @>>>
\HG \otimes_{\Lambda_{\Bbb Q_p}} H^1_{\Iw} (\Bbb Q_p,V(m))
$$
satisfying the following properties:

1) For any $n\geqslant 0$ the diagram
$$
\CD  \Cal D(V(m)) @>^
{\Exp^{\ep}_{V(m),h}}>>  
 \scr H (\Gamma) \otimes_{\Lambda_{\Bbb Q_p}} H^1_{\Iw} (\Bbb Q_p,V(m))\\
@V^{ \Xi^{\ep}_{V(m),n}}VV
@V^{\text{\rm pr}_{V(m),n}}VV
\\
H^1(K_n, C^{\bullet}_{\text{cris}}(\Ddagrig (V(m))) @>^{(h-1)!\exp_{V(m),K_n}}>>H^1(K_n,V(m))
\endCD
$$
commutes.

ii) Let $e_1= \ep^{-1}\otimes t$ denote the canonical generator of $\Dc (\Bbb Q_p(-1)).$ Then 
$$
\Exp_{V(m+1),h+1}^\ep=-\Tw^{\ep}_{V(m),1}\circ \Exp_{V(m),h}^{\ep}\circ (\partial \otimes e_1).
$$

iii) One has
$$
\Exp_{V(m),h+1}^\ep=\ell_h \Exp_{V(m),h}^\ep 
$$
where $\ell_m=m-\frac{\log (\gamma_1)}{\log \chi (\gamma_1)}.$
\endproclaim
\demo{Proof} This theorem was first proved  in \cite{PR2} for crystalline representations. Other proofs can be found in \cite{Cz1}, \cite{KKT}, \cite{Ben1} and \cite{Ber2}. Note that in \cite{Cz1} and \cite{KKT} $V$ is not assumed to be crystalline.
We also remark that in \cite{PR2} $\Exp_{V,h}^\ep (\alpha)$ was defined 
only for $\alpha$ such that $\partial^m\alpha (0)\in (1-p^m\Ph)\,\Dc (V)$ for all $m\in \Bbb Z.$   This condition is not necessary (see \cite{PR4} or \cite{Ben1, Section 5.1}).

 We recall now the construction of $\Exp_{V,h}^\ep$ in terms of $(\Ph,\Gamma)$-modules  found 
by Berger. This construction will be used in the proof of Proposition 2.2.2 below. 
Again in \cite{Ber2}, Berger assumes that $V$ is crystalline, but  his arguments work in the potentially
semistable case. We refer to \cite{Pt} for more detail. 
The action of $\scr H (\Gamma)$ on
$\bD^{\dag}(V)^{\psi=1}$ induces an isomorphism
$
\scr H (\Gamma) \otimes_{\La_{\Bbb Q_p}} \bD^{\dag}(V)^{\psi=1}
\rightarrow \Ddagrig (V)^{\psi=1}
$
(see  \cite{Pt, Section 6.4}).
Composing this map with the canonical isomorphism
$\Hi^1(\Bbb Q_p,V)\simeq \bD^{\dag}(V)^{\psi=1}$ we obtain an isomorphism
$
\scr H (\Gamma) \otimes_{\La_{\Bbb Q_p}} \Hi^1(\Bbb Q_p,V) \rightarrow \Ddagrig
(V)^{\psi=1}.
$
It is not difficult to check that $\ell_m$ acts on $\CR_L$ as $m- t\partial $ and an easy induction shows  that
$\underset{k=0}\to{\overset{h-1}\to \prod} \ell_k =(-1)^h t^h\dop^h.$
Let $h\geqslant 1$  be such that $\F^{-h}\Dd (V)=\Dd (V).$ 
To simplify the formulation, assume that $\Dc(V)^{\Ph=1}=0$. 
For any  $\alpha \in \Cal D(V)$ the equation
$$
(\Ph-1)\,F=\alpha -\underset{m=1}\to{\overset h \to \sum} \frac{\partial^m \alpha (0)}{m!}\, t^m 
$$
has a solution in $\CR_L^+\otimes \Dc (V)$ and we define 
$$
\Omega_{V,h}^\ep(\alpha) \,=\,\frac{\log \chi (\gamma_1)}{p}\,\ell_{h-1}\ell_{h-2} \cdots \ell_0 (F (X)).
$$
It is easy to see that $\Omega_{V,h}^\ep(\alpha) \in
\Drig^+(V)^{\psi=1}$ and in \cite{Ber2, Theorem II.13} Berger shows that
$\Omega_{V,h}^\ep(\alpha)$ coincides with $\Exp_{V,h}^\ep (\alpha).$
\enddemo
\flushpar
\,
\flushpar
{\bf 1.3.3. The logarithmic maps.} The Iwasawa algebra $\Lambda$ is equipped
with an involution $\iota \,:\,\Lambda @>>>\Lambda $ defined 
by $\iota (\tau)=\tau^{-1}$, $\tau \in \Gamma.$ If $M$ is a $\Lambda$-module
we set $M^{\iota}=\Lambda \otimes_{\iota}M$ and denote by $m\mapsto m^{\iota}$
the canonical bijection of $M$ onto $M^{\iota}.$ Thus $\lambda \, m^{\iota}=(\iota (\lambda)\,m)^{\iota}$ 
for all $\lambda \in \Lambda$, $m\in M.$ Let $T$ be a $O_L$-lattice of $V$ stable under the action
of $G_{\Bbb Q_p}.$ The cohomological pairings  
$$
(\,\,,\,\,)_{T,n}\,:\,H^1(K_n,T)\times H^1(K_n,T^*(1))@>>>O_L
$$
give rise to a $\Lambda$-bilinear pairing
$$
\bigl <\,\,\,,\,\,\bigr >_T\,:\,H^1_{\Iw}(\Bbb Q_p,T)\times H^1_{\Iw}(\Bbb Q_p,T^*(1))^{\iota} @>>>\Lambda 
$$
defined by 
$$
\bigl <x,y^{\iota}\bigr >_T\equiv \underset{\tau \in \Gamma/\Gamma_n}\to \sum (\tau^{-1}x_n,y_n)_{T,n}\tau \mod{(\gn-1)},
\qquad n\geqslant 1
$$
(see \cite{PR2, Section 3.6.1}). By linearity we extend this pairing to
$$
\bigl <\,\,\,,\,\,\bigr >_V\,:\,\HG\otimes_{\Lambda} H^1_{\Iw}(\Bbb Q_p,T)\times 
\HG \otimes_{\Lambda} H^1_{\Iw}(\Bbb Q_p,T^*(1))^{\iota} @>>>\HG .
$$
For any $\eta \in \Dc (V^*(1))$ the element $\widetilde \eta=\eta\otimes (1+X)$ lies in $\Cal D(V^*(1))$
and we define a map
$$
\Log^\ep_{V,1-h,\eta} \,\,:\,\, H^1_{\Iw}(\Bbb Q_p,V)@>>> \HG
$$
by
$$\Log^\ep_{V,1-h,\eta}(x)=\bigl <x,\Exp_{V^*(1),h}^{\ep^{-1}}(\tilde\eta)^{\iota} \bigr >_V.
$$
\proclaim{Lemma 1.3.4} For any $j\in \Bbb Z$ one has
$$
\Log^{\ep}_{V(-1),-h,\eta \otimes e_1}(\Tw_{V,-1}^\ep(x))\,=\,\Tw_1 \left (\Log^{\ep}_{V,1-h,\eta}(x)\right ).
$$
\endproclaim
\demo{Proof} A short computation shows that $\bigl <\Tw_{V,j}^{\ep}(x),\Tw_{V^*(1),-j}^\ep (y)\bigr >_{V(j)}=
\Tw_{-j}\bigl <x,y\bigr >_{V}.$ Taking into account that $\Tw_{V^*(1),1}^{\ep^{-1}}=-\Tw_{V^*(1),1}^{\ep}$ we have 
$$
\aligned
&\Log^{\ep}_{V(-1),-h,\eta \otimes e_1}(\Tw_{V,-1}^{\ep}(x))=
\bigl <\Tw_{V,-1}^{\ep}(x),  \Exp^{\ep^{-1}}_{V^*(2),h+1}(\widetilde{\eta\otimes e_1})^{\iota}\bigr >_{V(-1)}=\\
&\bigl <\Tw_{V,-1}^{\ep}(x), -\Tw_{V^*(1),1}^{\ep^{-1}}\bigl ( \Exp^{\ep^{-1}}_{V^*(1),h}(\widetilde \eta)\bigr )^{\iota}\bigr >_{V(-1)}=
\bigl <\Tw_{V,-1}^{\ep}(x), \Tw_{V^*(1),1}^{\ep} \bigl (\Exp^{\ep^{-1}}_{V^*(1),h}(\widetilde \eta)\bigr )^{\iota}\bigr >_{V(-1)}=\\
&\Tw_1 \bigl <x,  \Exp^{\ep^{-1}}_{V^*(1),h}(\widetilde \eta)^{\iota}\bigr >_{V}=
\Tw_1 \left (\Log^{\ep}_{V,1-h,\eta}(x)\right )
\endaligned
$$
and the lemma is proved.
\enddemo
\flushpar
{\bf 1.4. $p$-adic distributions} (see \cite{Cz6, Chapitre II}, \cite{PR2, Sections 1.1-1.2}).
Let $\Cal D (\Bbb Z^*_p,L)$ be the space of distributions on $\Bbb Z^*_p$ with values in a finite extensions $L$ of $\Bbb Q_p.$
To each  $\mu \in \Cal D(\Bbb Z_p^*,L)$ one can associate its Amice transform $\scr A_{\mu}(X)\in L[[X]]$ by
$$
\scr A_{\mu}(X)=\int_{\Bbb Z_p^*} (1+X)^x \mu (x)= \dsize \sum_{n=0}^\infty 
\left ( \int_{\Bbb Z_p^*}\binom{x}{n} \mu (x) \right )\,X^n.
$$
The map $\mu \mapsto \scr A_\mu (X)$ establishes an isomorphism between $\Cal D(\Bbb Z_p^*,L)$ and 
$(\CR^+_L)^{\psi=0}.$ We will denote by $\text{\rm \bf M}(\mu)$ the unique element of $\HG$ such that
$$
\text{\rm \bf M}(\mu)\,(1+X)\,=\,\Cal A_{\mu}(X).
$$ 
For each $m\in \Bbb Z$ the character $\chi^m\,:\,\Gamma @>>>\Bbb Z_p^*$ can be extended to a unique continuous
$L$-linear map  $\chi^m\,:\,\HG @>>>L^*$. If $h=\underset{i=0}\to{\overset{p-2}\to \sum} \delta_i h_i(\gamma_1-1)$,
then $\chi^m (h)=h_i(\chi^m(\gamma_1)-1)$ with $i\equiv m\pmod{(p-1)}.$ An easy computation
shows that
$$
\int_{\Bbb Z_p^*}x^m \mu (x)\,=\,\partial^m \scr A_{\mu}(0)\,=\,\chi^m(\text{\rm \bf M}(\mu)).
$$
If $x\in \Bbb Z_p^*$ we set $\left <x\right >=\omega^{-1}(x)\,x$ where $\omega$ denotes the Teichm\"uller character.
 To any $\mu \in \Cal D(\Bbb Z_p^*,L)$ we associate $p$-adic functions
$$
L_p(\mu,\omega^i,s)=\int_{\Bbb Z_p^*} \omega^i(x)\bigl <x\bigr >^s \mu (x), \qquad 0\leqslant i\leqslant p-2.
$$
Write $\text{\bf M}(\mu)=\underset{i=0}\to{\overset{p-2}\to\sum} \delta_i h_i(\gamma_1-1).$ Then
$$
L_p(\mu,\omega^i,s)=h_i \left ( \chi (\gamma_1)^{s}-1\right ). \tag{8}
$$
To prove this formula it is enough to compare the values of the both sides at the integers
\linebreak 
$s\equiv i\pmod{(p-1)}.$

We say that $\mu$ is of order $r>0$ if its Amice transform $\scr A_{\mu}(X)=\dsize \sum_{n=1}^\infty a_nX^n$ 
is of order $r$ i.e. if the sequence $|a_n|_p/n^r$ is bounded above. A distribution of order $r$ is completely 
determined by the values of the integrals
$$
\int_{\Bbb Z_p^*} \zeta_{p^n}^x x^i \mu (x),\qquad n\in \Bbb N, \quad 0\leqslant i\leqslant [r]
$$
where $[r]$ is the largest integer no greater then $r$.  

Set  $\hat {\Bbb Z}^{(p)}= \Bbb Z_p^*\times \underset{l\ne p}\to \prod \Bbb Z_{l}.$ 
A locally analytic function  on $\underset{l\ne p}\to\prod \Bbb Z_l$ is locally constant and we say
that a distribution $\mu$ on  $\hat {\Bbb Z}^{(p)}$ is of order $r$ if for any locally constant 
function $g(y)$ on $\underset{l\ne p}\to\prod \Bbb Z_l$ the linear map 
$f\mapsto \dsize\int_{\hat{\Bbb Z}^{(p)}} f(x) g(y) \mu (x,y)$ is a distribution of order $r$ on $\Bbb Z_p^*.$

\head
{\bf \S2. The $\ell$-invariant}
\endhead 

\flushpar
{\bf 2.1. The $\ell$-invariant.}  
\newline
{\bf 2.1.1. Definition of the $\ell$-invariant.} In this section we review and generalise slightly the definition of the $\ell$-invariant proposed in our previous article \cite{Ben2} in order to
cover the case of potentially crystalline reduction of modular forms. Let $S$ be a finite set of primes and $\Bbb Q^{(S)}/\Bbb Q$ be the 
maximal Galois extension of $\Bbb Q$ unramified outside $S\cup\{\infty\}.$ Fix a finite extension $L/\Bbb Q_p.$
Let $V$ be an $L$-adic  representation of $G_S$ i.e. a finite dimensional $L$-vector space equipped with a continuous
linear action of $G_S$.  We  write $H^*_S(\Bbb Q,V)$ for the continuous cohomology of $G_S$ with coefficients in $V$. 
We will always assume that the restriction of $V$ on the decomposition group at $p$
is potentially semistable.
For all  primes $l\ne p$ (resp. for $l= p$) Greenberg \cite{Gre1} (resp. Bloch and Kato \cite{BK}) defined a subgroup $H^1_f(\Bbb Q_l,V)$ of $H^1 (\Bbb Q_l,V)$ by
$$
H^1_f(\Bbb Q_l,V)=\cases \ker (H^1(\Bbb Q_l,V)@>>>H^1(\Bbb Q_l^{\text{\rm ur}},V)) &\text{if $l\ne p$,}\\
 \ker (H^1(\Bbb Q_p,V)@>>>H^1(\Bbb Q_p,V\otimes \Bc)) &\text{if $l=p$}
\endcases
$$
where $\Bc$ is the ring of crystalline periods \cite{Fo2}.
The Selmer group of $V$ is defined as 
$$
H^1_f\left (\Bbb Q,V \right )=\ker \left (H^1_S(\Bbb Q, V) @>>> \underset{l\in S}\to \bigoplus \frac{H^1(\Bbb Q_l,V)}{H^1_f(\Bbb Q_l,V)}\right ).
$$ 
We also define
$$
H^1_{f,\{p\}}\left (\Bbb Q,V \right )=\ker \left (H^1_S(\Bbb Q, V) @>>> \underset{l\in S-\{p\}}\to \bigoplus \frac{H^1(\Bbb Q_l,V)}{H^1_f(\Bbb Q_l,V)}\right ). 
$$

Note that these definitions do not depend on the choice of $S$.  From now until the end of this \S   we assume that $V$
satisfies the following conditions

{\bf 1)}  $H^1_f(\Bbb Q,V)=H^1_f(\Bbb Q,V^*(1))=0.$ 

 {\bf 2)} The action of $\Ph$ on $\Dst (V)$ is semisimple at $1$.

{\bf 3)}  ${\text{\rm dim}}_{L} t_V(\Bbb Q_p)=1.$

 We remark that the last condition can be relaxed but it simplifies the formulation of Proposition 2.2.4 below and holds for the situations considered in \S\S3-4.

The condition {\bf 1)} together with the Poitou-Tate exact sequence (see \cite{FP, Proposition 2.2.1})
$$
\cdots @>>>H^1_f(\Bbb Q,V)@>>>H^1_S(\Bbb Q,V)\simeq \underset{l\in S}\to \bigoplus \frac{H^1(\Bbb Q_l,V)}{H^1_f(\Bbb Q_l,V)}
@>>>H^1_f(\Bbb Q,V_f^*(1))^*@>>>\cdots
$$
gives an isomorphism
$$
H^1_S(\Bbb Q,V)\simeq \underset{l\in S}\to \bigoplus \frac{H^1(\Bbb Q_l,V)}{H^1_f(\Bbb Q_l,V)}. 
$$
In particular, we have 
$$
H^1_{f,\{p\}}(\Bbb Q,V)\simeq \frac{H^1(\Bbb Q_p,V)}{H^1_f(\Bbb Q_p,V)}. \tag{9}
$$

Let $D$ be a one-dimensional subspace of $\Dc (V)$ on which $\Ph$ acts as multiplication by $p^{-1}.$
Set $D_{\Bbb Q_p^{\text{\rm ur}}}=D\otimes_{\Bbb Q_p}\Bbb Q_p^{\text{\rm ur}}$ where $\Bbb Q_p^{\text{\rm ur}}$ denotes the
maximal unramified extension of $\Bbb Q_p.$ Using the weak admissibility of $\Dpst (V)$ it is easy to see that  $D$ is not contained in $\F^0\Dc (V)$ and therefore
$$\Dpst (V)=\F^0\Dpst (V)\oplus D_{\Bbb Q_p^{\text{\rm ur}}} \tag{10}
$$
as $\Bbb Q_p^{\text{\rm ur}}\otimes_{\Bbb Q_p}L$-modules.  Let $m$ denote the unique Hodge-Tate weight of $D$.
By Berger's theory \cite{Ber4} (see also \cite{BC, Section 2.4.2}), the intersection $\Ddagrig (V)\cap \bigl(D
\otimes_L \CR_L[1/t]\bigr )$ is a saturated $(\Ph,\Gamma)$-submodule of  $\Ddagrig (V)$ of rank $1$ which is 
isomorphic to $\CR_L(\delta)$ with $\delta (x)=|x|x^{m}.$ Thus we have an exact sequence of $(\Ph,\Gamma)$-modules
$$
0@>>> \CR_L(\delta)@>>>\Ddagrig (V) @>>> \bold D@>>>0 \tag{11}
$$ 
where $\bold D=\Ddagrig (V)/\CR_L(\delta).$ Passing to duals and taking the 
long exact cohomology sequence we obtain an exact sequence
$$
 H^1(\Bbb Q_p,V^*(1))@>>>H^1(\CR_L(\chi\delta^{-1}))@>>>H^2(\bold D^*(\chi)). \tag{12}
$$

\proclaim{Proposition 2.1.2} Assume that one of the following conditions holds

a) $D$ is not contained in the image of the monodromy operator $N\,:\,\Dpst (V)@>>>\Dpst (V)$
and $\Dc (V)^{\Ph=1}=0.$

b) $D$ is contained in the image of $N$ and $N^{-1}(D)\cap \Dst (V)^{\Ph=1}$ is a one-dimensional
$L$-vector space.
\flushpar
Then the composition
$$
\varkappa\,:\,H^1_{f,\{p\}}(\Bbb Q, V^*(1))
@>>> H^1(\CR_L(\chi\delta^{-1}))
$$
 of the localisation map 
$H^1_{f,\{p\}}(\Bbb Q, V^*(1))@>>>H^1(\Bbb Q_p, V^*(1))$ with 
 $H^1(\Bbb Q_p,V^*(1))@>>>H^1(\CR_L(\chi\delta^{-1}))$ is injective.
Moreover, $\text{\rm Im}(\varkappa)$ is a one-dimensional $L$-vector
space such that 
$$
\text{\rm Im}(\varkappa) \cap H^1_f( \CR_L(\chi\delta^{-1}))=\{0\}.
$$
\endproclaim
\demo{Proof} We consider the two cases separately.

a) First assume that $D$ satisfies a). Applying the functor $\CDpst$ to (11) we obtain
an exact sequence
$$
0@>>>D_{\Bbb Q_p^{\text{\rm ur}}} @>>>\Dpst (V)@>f>>\CDpst (\bold D)@>>>0. \tag{13}
$$
From (10) it follows that $\F^0 \CDpst (\bold D)=\CDpst (\bold D)$ and by \cite{Ben2, Proposition 1.4.4} 
$$H^0(\bold D)\simeq \CDpst (\bold D)^{\Ph=1,N=0, G_{\Bbb Q_p}} =\CDcris (\bold D)^{\Ph=1}.
$$
 Applying the snake lemma to (13) we obtain an isomorphism
of $G_{\Bbb Q_p}$-modules $\Dpst (V)^{\Ph=1} \simeq \CDpst (\bold D)^{\Ph=1}.$ Thus 
$\Dst (V)^{\Ph=1} \simeq \CDst (\bold D)^{\Ph=1}.$
Let $x\in \CDcris (\bold D)^{\Ph=1}.$
There exists a unique $y\in \Dst (V)^{\Ph=1}$ such that
$f(y)=x.$ Since $f(N(y))=N(x)=0$, one has $N(y)\in D$ and 
by a) $N(y)=0$ i.e. $y\in \Dc (V)^{\Ph=1}.$  Since  $\Dc (V)^{\Ph=1}=0$ by assumption a), we proved that 
$H^0(\bold D)=0.$

Now $H^2(\bold D^*(\chi))=0$ by Poincar\'e duality
and from the exact sequence (12) we obtain that the map
$H^1(\Bbb Q_p,V^*(1))@>>>H^1(\CR_L(\chi\delta^{-1}))$ is surjective. Since
$\chi\delta^{-1}(x)=x^{1-m},$ the cohomology $H^1(\CR_L(\chi\delta^{-1}))$ decomposes
into the direct sum of one dimensional subspaces
$$
H^1(\CR_L(\chi\delta^{-1}))\simeq H^1_f(\CR_L(\chi\delta^{-1}))\oplus H^1_c(\CR_L(\chi\delta^{-1})).
$$
The image of $H^1_f(\Bbb Q_p,V^*(1))$ in $H^1(\CR_L(\chi\delta^{-1}))$ is contained 
in  $H^1_f(\CR_L(\chi\delta^{-1}))$ and we have a surjective map
$$
\frac{H^1(\Bbb Q_p,V^*(1))}{H^1_f(\Bbb Q_p,V^*(1))}@>>>
\frac{H^1(\CR_L(\chi\delta^{-1}))}
{H^1_f(\CR_L(\chi\delta^{-1}))}. \tag{14}
$$
From $\Dc (V)^{\Ph=1}=0$ it follows that
$H^0(\Bbb Q_p,V)=0$ and 
$$
\text{\rm dim}_L \left (H^1_f(\Bbb Q_p,V)\right )=\text{\rm dim}_L \left (t_V(L)
\right )+
\text{\rm dim}_L \left (H^0(\Bbb Q_p,V)\right )= 1.
$$ 
Therefore $H^1(\Bbb Q_p,V^*(1))/H^1_f(\Bbb Q_p,V^*(1))$ is one-dimensional 
and the map (14) is an isomorphism. The Proposition follows now from
this fact and from the isomorphism (9) for the cohomology with coefficients in $V^*(1)$ 
instead $V$.

b) Now assume  that $D$ satisfies b). We follow the approach of \cite{Ben2, Sections 2.1 and
2.2} with some modifications ( see especially the proofs   of Proposition 2.1.7 and
Lemma 2.1.8 of op. cit.). The Proposition will be proved in several steps.

b1) Consider the filtration on $\Dpst (V)$ given by
$$
D_i=\cases 0 &\text{if $i=-1$}\\
D_{\Bbb Q_p^{\text{\rm ur}}} &\text{if $i=0$}\\
(D+N^{-1}(D)\cap \Dst (V)^{\Ph=1}) _{\Bbb Q_p^{\text{\rm ur}}} &\text{if $i=1$}\\
\Dpst (V) &\text{if $i=2$.}
\endcases
$$
By \cite{Ber4} this filtration induces a unique filtration on $\Ddagrig (V)$
$$
\{0\}= F_{-1}\Ddagrig (V) \subset F_0\Ddagrig (V)\subset F_1\Ddagrig (V)\subset F_2\Ddagrig (V)=\Ddagrig (V)
$$ 
such that $\CDpst (F_i\Ddagrig (V))=D_i.$ Note that $F_1\Ddagrig (V)$ is a semistable $(\Ph,\Gamma)$-submodule
of $\Ddagrig (V)$. To simplify notation set $M_0=F_0\Ddagrig (V)$, $M=F_1\Ddagrig (V)$ and
$M_1=\text{\rm gr}_1\Ddagrig (V).$ We remark that $M_0\simeq \CR_L(\delta)$ and since 
$\F^0(D_1/D_0)=D_1/D_0$ and $(D_0/D_1)^{\Ph=1}=D_0/D_1$ Proposition 1.5.9 of \cite{Ben2} implies that
$M_1\simeq \CR_L(x^{-k})$ for some $k\geqslant 0.$ By the assumption b) the monodromy operator
$N$ acts non trivially on $D_1$ and therefore we have a non crystalline extension
$$
0@>>>\CR_L(\delta)@>>>M@>>>\CR_L(x^{-k})@>>>0
$$
which is a particular case of the exact sequence from \cite{Ben2, Proposition 2.1.7}. Passing to duals and taking
the long cohomology sequence we obtain a diagram
$$
\xymatrix{
H^1(\CR_L(\vert x\vert x^{k+1})) \ar[r]^{f_1} & H^1(M^*(\chi)) \ar[r]^{g_1} & H^1(\CR_L(\chi\delta^{-1})) \ar[r]^{\delta_1} & H^2(\CR_L(\vert x\vert x^{k+1})) \ar[r] & 0\,.\\
& H^1_{f,\{p\}}(\Bbb Q,V^*(1))\ar[u]^{\eta} \ar[ur]^{ \varkappa} &  & &
} \tag{15}
$$ 

b2) The quotient $\bD=\text{\rm gr}_2\Ddagrig (V)$ is a potentially semistable $(\Ph,\Gamma)$-module with Hodge-Tate
weights $\geqslant 0$. Thus $H^0\left (\bD \right )\simeq
( \Dpst (V)/D_1)^{G_{\Bbb Q_p},\Ph=1,N=0}.$ Let $\bar x=x+D_1\in ( \Dpst (V)/D_1)^{G_{\Bbb Q_p},\Ph=1,N=0}.$
For each $g\in G_{\Bbb Q_p}$ we can write $g(x)=x+d$ for some $d\in D_1.$
Since the inertia subgroup $I_p\subset G_{\Bbb Q_p}$ acts on $\Dpst (V)$ through a finite quotient and 
since the restriction of this action on  $D_1$ is trivial, we obtain that 
$x\in \Dpst (V)^{I_p}=\Dst (V)\otimes_{\Bbb Q_p}\Bbb Q_p^{ur}.$ Thus
$$
\bar x\in \left (\left (\frac{\Dst(V)}{D+N^{-1}(D)\cap \Dst (V)^{\Ph=1}}\right )
\otimes{\Bbb Q_p^{ur}}\right )^{G_{\Bbb Q_p},\Ph=1,N=0}=
\left (\frac{\Dst(V)}{D+N^{-1}(D)\cap \Dst (V)^{\Ph=1}}\right )^{\Ph=1,N=0}.
$$ 
Since $\Ph$ is semisimple at $1$, we can assume that $x\in \Dst (V)^{\Ph=1}.$
Then $N(x)\in D$ and therefore $x\in N^{-1}(D)\cap \Dst (V)^{\Ph=1}.$ 
This shows that $\bar x=0$ and we proved that $H^0\left (\bD\right )=0.$
By Poincar\'e duality we obtain immediately that $H^2(\bD^*(\chi))=0.$
Now the Euler characteristic formula together with \cite{Ben2, Corollary 1.4.5} give
$$
\dim_L H^1 (\bD^* (\chi) ))=
\text{\rm rg} \left (\bD^*(\chi)\right )+\dim_L H^0\left ( \bD^*(\chi)\right )=\dim_L H^1_f \left (\bD^* (\chi)\right )
$$
and therefore  
$$
H^1_f \left (\bD^* (\chi)\right )=
H^1 \left (\bD^* (\chi)\right ). \tag{16}
$$

b3) Consider the exact sequence 
$$
0@>>>\bD^*(\chi)@>>>\Ddagrig (V^*(1))@>>>M^*(\chi)@>>>0.
$$
Since $H^0(M^*(\chi))=0$, this sequence together with the isomorphism (16)  give an  exact sequence
$$
0@>>>H^1_f\left (\bD^*(\chi) \right)    @>>>H^1(\Bbb Q_p,V^*(1))@>>>H^1(M^*(\chi))@>>>0. \tag{17}
$$
On the other hand, by \cite{Ben2, Corollary 1.4.6} the sequence
$$
0@>>>H^1_f\left (\bD^*(\chi) \right)    @>>>H^1_f(\Bbb Q_p,V^*(1))@>>>H^1_f(M^*(\chi))@>>>0 \tag{18}
$$
is also exact. 

b4) We come back to the diagram (15). The sequences (17) and (18) together with the isomorphism (9) show that
the map $\eta$ is injective. 
By \cite{Ben2, Lemma 2.1.8} one has $\ker (g_1)=H^1_f(M^*(\chi))$ and 
$H^1(\CR_L(\chi\delta^{-1}))=H^1_f(\CR_L(\chi\delta^{-1})) \oplus   \text{\rm Im}(g_1).$
Since 
$$
H^1(\Bbb Q_p,V^*(1))/H^1_f(\Bbb Q_p,V^*(1))\simeq H^1(M^*(\chi))/H^1_f(M^*(\chi))
$$
we obtain that  $\ker (\varkappa)=\text{\rm Im}(\eta)\cap \ker (g_1)=0$ and that  
$\text{\rm Im}(\varkappa) \cap H^1_f(\CR_L(\chi\delta^{-1}))=0.$
The Proposition is proved. 
\enddemo

\proclaim{Definition}  The $\ell$-invariant associated to $V$ and $D$
is the unique element  $\ell (V,D)\in L$ such that
$$
\text{\rm Im} (\varkappa) = L(\bold y_{m-1}+\ell (V,D)\,\bold x_{m-1}).
$$
 Here $\{\bold x_{m-1},\bold y_{m-1}\}$
is the canonical basis of $H^1(\CR_L(\chi\delta^{-1}))$ constructed in 1.2.5.
\endproclaim

{\bf Remarks 2.1.3.} 1) If $V$ is semistable at $p$ this definition agrees with
the definition of $\ell (V,D)$ proposed in \cite{Ben2, Sections 2.2.2 and 2.3.3}.

2) The assumptions a) and b) imply both that $H^0(\Bbb Q_p,V)=0.$

3) One can express $\ell (V,D)$ directly in terms of $V$ and $D$. 
Assume that $D$ satisfies the condition a) of Proposition 2.1.2.
Since  $H^0(\bold D)=0$ the sequence (11) shows that
$H^1(\CR_L(\delta))$ injects into $H^1(\Ddagrig (V))\simeq H^1(\Bbb Q_p,V).$
Moreover, from 
$
\text{\rm dim}_LH^1_f(\Bbb Q_p,V)=1
$
and the fact that $\dim_LH^1_f(\CR_L(\delta))=1$  it follows that 
$H^1_f(\Bbb Q_p,V)\simeq H^1_f(\CR_L(\delta)).$ Let $H^1_D(V)$ denote the inverse image of
$H^1(\CR_L(\delta))/H^1_f(\CR_L(\delta))$ under the isomorphism (9).  Then
$$
H^1_D(V)\simeq \frac{H^1(\CR_L(\delta))}{H^1_f(\CR_L(\delta))} 
$$
and  the localisation map
$H^1_S(V)@>>>H^1(\Bbb Q_p,V)$ induces an injection $H^1_D(V)@>>>H^1(\CR_L(\delta)).$  
Using the decomposition (6) we define $\scr L(V,D)$ as the unique element of $L$ such that
$$
\text{\rm Im} (H^1_D(V)@>>>H^1(\CR_L(\delta)))=L(\boldsymbol \beta_m+\scr L(V,D) \boldsymbol \alpha_m)
$$
where $\{\boldsymbol \alpha_m,\boldsymbol \beta_m\}$ denotes the canonical basis of $H^1(\CR_L(\delta))$.
Then
$$
\ell (V,D)=-\scr L(V,D) \tag{19}
$$
(see \cite{Ben2, Proposition 2.2.7}). Note that in op. cit. $V$ is assumed to be semistable,
but in the potentially semistable case the proof is exactly the same.  

A similar duality formula can be proved in the case b) too, but it will not be used in this paper.
We refer to \cite{Ben2, Section 2.2.3} for more detail.

4) The diagram (15) shows that in the case b)  the image of $H^1_{f,\{p\}}(\Bbb Q,V^*(1))$ in 
$H^1(\CR_L(\chi\delta^{-1}))$ coincides with  $\text{\rm Im}(g_1)$ and therefore that $\ell (V,D)$ 
depends only on the local properties of $V$ at $p$. On the other hand, in the case a)
the $\ell$ invariant is global and contains information about the localisation map
$H^1_{f,\{p\}}(\Bbb Q,V^*(1))@>>>H^1(\CR_L(\chi\delta^{-1})).$
\newline
\newline
{\bf 2.2. Relation to the large exponential map.}
\newline
{\bf 2.2.1. Derivative of the large exponential map.} In this section we interpret  $\ell (V,D)$ in terms
of the Bockstein homomorphism associated to the large exponential map.
This interpretation is crucial for the  proof of the main theorem of this paper.
We keep the notations and conventions of Section 2.1.
Recall (see Section 1.3.2) that  $H^1(\Bbb Q_p,\scr H(\Gamma)\otimes_{\Bbb Q_p}V)= \scr H(\Gamma)\otimes_{\Lambda_{\Bbb Q_p}} H_{\Iw}^1(\Bbb Q_p,V)$
injects into  $\Ddagrig (V).$  Set 
$$H^1_{\delta}(\Bbb Q_p,\scr H(\Gamma)\otimes_{\Bbb Q_p}V)=  \CR_L(\delta)
\cap H^1(\Bbb Q_p,\scr H(\Gamma)\otimes_{\Bbb Q_p}V). $$ 
The projection map induces a commutative diagram
$$
\CD
H^1_{\delta}(\Bbb Q_p,\scr H(\Gamma)\otimes_{\Bbb Q_p}V) @>>> H^1(\Bbb Q_p,\scr H(\Gamma)\otimes_{\Bbb Q_p}V)\\
@VVV @V\text{\rm pr}_{0}VV\\
H^1(\CR_L(\delta)) @>>> H^1(\Bbb Q_p,V)
\endCD
$$
where the bottom arrow is an injection.  We fix a generator $\gamma \in \Gamma$
and an integer $h\geqslant 1$ such that  $\F^{-h}\Dd (V)=\Dd (V).$

\proclaim{Proposition 2.2.2}  Assume that $D$ is a one-dimensional subspace of $\Dc (V)$ on which $\Ph$ acts
as multiplication by $p^{-1}.$  
For any  $a \in D$ let $x  \in \Cal D(V)$
be such that $x (0)=a .$ Then

i) There exists a unique $F \in  H^1_{\delta}(\Bbb Q_p,\scr H(\Gamma)\otimes V)$ such that
$$
(\gamma-1)\,F =\Exp_{V,h}^\ep (x ).
$$

ii) The composition map
$$
\align
&\delta_{D}\,:\,D@>>> H^1_{\delta}(\Bbb Q_p,\scr
H(\Gamma)\otimes V)@>>>H^1(\CR_L(\delta))\\
&\delta_{D} (a)=\text{\rm pr}_{0}(F)
\endalign
$$
is well defined and is explicitly given
 by  the following formula 
$$
\delta_{D}(a)\,=\,\Gamma \left (h\right )\,\left (1-\frac{1}{p}\right )^{
 -1}(\log \chi
(\gamma))^{-1}\, i_{c} (a).
$$
\endproclaim
\demo{Proof} 1) Since  in both cases, a) and b) $\Dc (V)^{\Ph=1}=0,$ the operator $1-\Ph$ is invertible on
$\Dc (V)$ and we have a diagram
$$
\xymatrix{ \Cal D(V) \ar[rr]^{\Exp_{V,h}^\ep}
\ar[d]^{\Xi_{V,0}^\ep} & & H^1(\Bbb Q_p, \scr H(\Gamma)\otimes V)
\ar[d]^{\text{\rm pr}_{0}}\\
t_{V}(\Bbb Q_p)\otimes \Dc (V) \ar[rr]^{(h-1)!\exp_{V}} & & H^1(\Bbb Q_p,V).}
$$
where $\Xi_{V,0}^\ep
(f)=\dsize \left (\frac{1-p^{-1}\Ph^{-1}}{1-\Ph}\,f (0),0\right )$. 
If
$x \in D \otimes \CR_L^{\psi=0}$ then
$\Xi_{V,0}^\ep(x)=0$ and 
$
\text{\rm pr}_{0} \left (\Exp_{V,h}^\ep (x)\right )\,=\,0.
$
On the other hand, as $H^1_{\Iw}(\Bbb Q_p,V)$ is $\Lambda_{\Bbb Q_p}$-free,  the
 map 
$
\left (\scr H(\Gamma)\otimes_{\Lambda_{\Bbb Q_p}} \Hi^1(\Bbb Q_p,V)\right )_{\Gamma} @>>> H^1(\Bbb Q_p,V)
$
is injective and therefore  there exists a unique 
$F \in \scr H(\Gamma)\otimes_{\Lambda_{\Bbb Q_p}} \Hi^1 (\Bbb Q_p,V)$
such that 
$
\Exp_{V,h}^\ep (x) \,=\,(\gamma-1)\,F.
$
Let  $y\in D\otimes \CR_L^{\psi=0}$ be  another element such that $y(0)=a$ and let
$\Exp_{V,h}^\ep (y)=(\gamma-1)\,G.$  
Since  $\CR_L^{\psi=0}=\HG \,(1+X)$ we have 
$y=x+(\gamma-1)g$ for some $g\in D\otimes \CR_L^{\psi=0}.$ As $\Exp_{V,h}^\ep (g)=0$, we obtain immediately
that $\text{\rm pr}_{0}(G)=\text{\rm pr}_{0}(F)$ and we proved that the map $\delta_{D}$ is well defined.  

2) Take $a\in D$ and set
$$
x= a \otimes \ell \left (\frac{(1+X)^{\chi (\gamma)}-1}{X}\right ),
$$
where  
$
\ell (u)\,=\,\dsize\frac{1}{p}\log \left(\frac{u^p}{\Ph (u)}\right ). 
$
 An easy computation shows that
$$
\sum_{\zeta^p=1} \ell \left ( \frac{\zeta^{\chi (\gamma)}(1+X)^{\chi (\gamma)}-1}{\zeta (1+X)-1} \right )\,=\,0.
$$
Thus $x \in D \otimes O_L[[X]]^{\psi=0}.$
Write $x$ in the form
$f=(1-\Ph)\,(\gamma-1)\,(a\otimes \log (X)).$ Then 
$$
\Omega_{V,h}^\ep(x)\,=\,(-1)^{h-1}\frac{\log \chi (\gamma_1)}{p}\,t^{h}\partial^{h}((\gamma-1)\,(a\log (X))\,=\,
\left (1-\frac{1}{p}\right ) \log \chi (\gamma)\,(\gamma-1)\,F
$$
where 
$$
F \,=\,(-1)^{h-1}t^{h} \partial^{h} (a\log (X))\,=\,(-1)^{h-1}a t^{h} \partial^{h-1}
\left (\frac{1+X}{X}\right ).
$$
This implies immediately that $F \in  H^1_{\delta} (\Bbb Q_p, \scr H (\Gamma)\otimes V).$ On the other hand, as 
$
D\,=\,\Cal D_{\text{\rm cris}}(\CR_L(\delta))
$
without lost of generality  we may assume that $a\,=\,t^{-m}e_{\delta}$ where $\delta(x)=\vert x\vert x^m.$
Then 
$$F= (-1)^{h-1}t^{h-m}\partial^{h} \log (X)\,e_\delta.
$$
One has 
 $\text{\rm pr}_{0}(F)=\cl (G,F)$ where 
$(1-\gamma)\,G\,=\,(1-\Ph)\,F$ (see Section 1.2.3) and  
by \cite{CC1, Lemma 1.5.1}
there exists a unique $b\in \CE^{\dag, \psi=0}_{L}$ 
such that $(1-\g)\,b=\ell (X)$. One has  
 $$
(1-\gamma)\,\left (t^{h-m}\partial^{h} be_{\delta}\right )\,=\,(1-\Ph)\,\left (t^{h-m}\partial^{h} \log (X)e_{\delta}\right )\,=\,(-1)^{h-1}
(1-\Ph)\,F.
$$
Thus $G= (-1)^{h-1}t^{h-m}\partial^{h} (b) e_{\delta}$ and 
$
\res \left (G\,t^{m-1}dt \right )=(-1)^{h-1}\res \left (t^{h-1}\partial^{h} (b) dt \right )\,e_{\delta}=0.
$
Next from the congruence
$
F \equiv (h-1)!\,t^{-m}e_{\delta}\,
\pmod{\Bbb Q_p[[X]]\,e_{\delta}}
$
it follows  that
$\res (F t^{m-1}dt)\,=\,(h-1)!\,e_{\delta}.$ Therefore by \cite{Ben2, Corollary 1.5.5}  we have 
$$
\left (1-\frac{1}{p} \right )\,(\log \chi (\g))\,\cl (G,F)\,=\,(h-1)!\,\cl (\boldsymbol\beta_m)\,=\,(h-1)!   \,i_{c}(a).\tag{20}
$$
 On the other hand 
$$
x (0)\,=\,\left.
a \otimes \ell \left (\frac{(1+X)^{\chi (\gamma)}-1}{X}\right )
 \right \vert_{X=0}\,=\,a\left (1-\frac{1}{p} \right )\,\log \chi (\gamma).\tag{21}
$$
The formulas (20) and (21) imply that
$$
\delta_{D}(a)=(h-1)!\,\left (1-\frac{1}{p}\right )^{
 -1}(\log \chi
(\gamma))^{-1}\, i_{c} (a).
$$
and the Proposition is proved.
\newline
\,
\enddemo
\,
\flushpar
{\bf 2.2.3. Derivative of the large logarithmic map and $\ell$-invariant.} Fix a non-zero element $d\in D$ and consider the large logarithmic map 
$$
\Log^{\ep}_{V^*(1),1-h,d}\,\,:\,\,H^1_{\Iw}(\Bbb Q_p,V^*(1))@>>>\HG
$$ 
 (see Section 1.3.3). Let 
$$
H^1_{\Iw,S}(\Bbb Q,T^*(1))=\underset{\text{\rm cor}}\to \varprojlim H^1_S(\Bbb Q (\zeta_{p^n}),T^*(1))
$$
denote the global Iwasawa cohomology with coefficients in $T^*(1)$ and let 
$$H^1_{\Iw,S}(\Bbb Q,V^*(1))=H^1_{\Iw,S}(\Bbb Q,T^*(1))\otimes_{\Bbb Z_p} \Bbb Q_p.$$
The main results of this paper will be directly deduced 
from the following statement.

\proclaim{Proposition 2.2.4} Assume that $D$ is a one-dimensional subspace of $\Dc (V)$ on which $\Ph$ acts
as multiplication by $p^{-1}$ and which satisfies one of the conditions a-b) of Proposition 2.1.2. 
Let $\bold z\in H^1_{\Iw,S}(\Bbb Q,V^*(1))$. Assume that 
$\bold z_0=\text{\rm pr}_0(\bold z)\in H^1_S(\Bbb Q,V^*(1))$
is non-zero and denote by $\mu_{\bold z}\in \Cal D(\Bbb Z_p^*,L)$   the distribution
defined by
$$
{\bold M}(\mu_{\bold z})=\Log^{\ep}_{V^*(1),1-h,d}(\bold z).
$$
Consider the $p$-adic function 
$$
L_p(\mu_{\bold z},s)=\int_{\Bbb Z_p^*}\bigl <x\bigr >^s\mu_{\bold z}(x).
$$
Then $L_p(\mu_{\bold z},0)=0$ and
$$
L'_p(\mu_{\bold z},0)\,=\,\ell (V,D)\,\, \Gamma \left (h\right ) \left (1-\frac{1}{p}\right )^{-1}\,
 \left [d, \exp_{V^*(1)}^*({\bold z}_0)\right ]_{V}
$$ 
where $[\,\,,\,]_{V}\,:\,\Dc (V)\times \Dc (V^*(1))@>>>L$ is the canonical duality.
\endproclaim
\demo{Proof} First note that by \cite{PR1, Section 2.1.7}
for $l\neq p$ one has  $H^1_{\Iw}(\Bbb Q_{l},V^*(1))\simeq H^0(\Bbb Q_l(\zeta_{p^\infty}),V^*(1))$
 and therefore  ${H_{\Iw}^1(\Bbb Q_l,V^*(1))}_{\Gamma}$ is contained in $H^1_f(\Bbb Q_l,V^*(1)).$ Thus 
${H_{\Iw,S}^1(\Bbb Q,V^*(1))}_{\Gamma}$ injects into 
\linebreak
$H^1_{f,\{p\}}(\Bbb Q,V^*(1))$ 
and $\bold z_0\in H^1_{f,\{p\}}(\Bbb Q,V^*(1)).$ Recall that we fixed a basis $d$  of the one-dimensional
$L$-vector space $D=\CDcris (\CR_L(\delta)).$ Let $d^*$ be the basis of 
$\CDcris (\CR_L(\chi\delta^{-1}))$ which is dual to $d.$ Let $\tilde{\bold z}_0$ denote the image of 
$\bold z_0$ under the projection map $H^1(\Ddagrig (V^*(1)))@>>>H^1(\CR_L(\chi\delta^{-1})).$
 Write $\tilde{\bold z}_0=a\, i_f(d^*)+b\,i_c(d^*).$ Then $\ell (V,D)=a/b$.
By Proposition 1.2.6 and (7) we have 
$$
\multline
\left [d, \exp_{V^*(1)}^*(\bold z_0)\right ]_{V}\,=\,-
\exp_{V} (d)\cup \,\bold z_0\,=\,-\exp_{\CR_L(\delta)}(d)\cup \,\,\tilde{\bold z}_0\,=\,\\
=-b\,\bigl (i_f(d)\cup i_c(d^*)\bigr )\,=\,
-b\,\left (\boldsymbol \alpha_{m}\cup \bold y_{m-1}\right )\,=\,b.
\endmultline \tag{22}
$$
Let $\bold M (\mu_{\bold z})\,=\,\underset{i=0}\to{\overset{p-2}\to \sum}\delta_i h_i(\gamma_1-1).$
Then   $L_p(\mu_{\bold z},s)=h_0 \left (\chi (\gamma_1)^s-1 \right )$ by (8).
From Proposition 2.2.2 it follows that there exists $F\in H^1_{\delta}(\Bbb Q_p,\scr H(\Gamma)\otimes V)$
such that
$\Exp^{\ep^{-1}}_{V,h}( d\otimes (1+X))= (\gamma-1)\,F$ and  
$$
\bold M(\mu_{\bold z})=\Log_{V^*(1),1-h,d}^{\ep}(\bold z)=
\left <\bold z\,,\,\Exp^{\ep^{-1}}_{V,h}( d\otimes (1+X)^{\iota}\right >_{V}\,=\,
(\gamma^{-1}-1)\,\left <\bold z, F^{\iota} \right >_{V}.
$$ 
Put $\left <\bold z, F^{\iota} \right >_{V_f}=\underset{i=0}\to{\overset{p-2}\to \sum} \delta_i H_i(\gamma_1-1).$ Then $L_p(\mu_{\bold z},s)=(\chi (\gamma)^{-s}-1)\,H_0(\chi (\gamma_1)^s-1).$ Since $\chi (\gamma_1)=\chi (\gamma)^{p-1}$
the last formula implies that $L_p(\mu_{\bold z},s)$ has a zero at $s=0$ and 
$$
L_p'(\mu_{\bold z},0)=-(\log \chi (\g))\,H_0(0). \tag{23}
$$
On the other hand, by Proposition 2.2.2 
$$
\multline 
H_0(0)=\bold z_0 \cup (\text{\rm pr}_0F)=\tilde{\bold z}_0\cup \delta_D (d)=\Gamma \left (h\right ) \left (1-\frac{1}{p}\right )^{-1} 
(\log \chi(\gamma))^{-1}\, \bigl (\tilde{\bold z}_0 \cup i_c(d)\bigr )\,=\,\\ 
=-\Gamma \left (h\right ) \left (1-\frac{1}{p}\right )^{-1} 
(\log \chi(\gamma))^{-1}a . 
\endmultline \tag{24}
$$
From (22), (23) and (24) we obtain that 
$$
L_p'(\mu_{\bold z},0)=\Gamma \left (h\right ) \left (1-\frac{1}{p}\right )^{-1} a\,=\,
\ell (V,D)\,\, \Gamma \left (h\right ) \left (1-\frac{1}{p}\right )^{-1}\,
\left [d_\alpha, \exp_{V^*(1)}^*(\bold z_0)\right ]_{V}
$$
and the Proposition is proved.
\enddemo

\head
{\bf \S3. Trivial zeros of Dirichlet $L$-functions}
\endhead

{\bf 3.1. Dirichlet $L$-functions.}  Let $\eta\,:\,(\Bbb Z/N\Bbb Z)^*@>>>\Bbb C^*$ be a Dirichlet character of conductor $N$. 
We fix a primitive $N$-th root of unity $\zeta_N$ and set
$\tau (\eta)=\underset{a\mod N}\to \sum \eta (a)\,\zeta_N^a$.
The Dirichlet $L$-function
$$
L(\eta,s)=\underset{n=1}\to{\overset\infty\to \sum}\frac{\eta (n)}{n^s},\qquad \text{\rm Re}(s)>1
$$
has a meromorphic continuation on the whole complex plane and satisfies the functional equation
$$
\left (\frac{N}{\pi}\right )^{s/2}\Gamma \left (\frac{s+\delta_\eta}{2}\right )\,L(\eta,s)\,=\,
W_\eta  \left (\frac{N}{\pi}\right )^{(1-s)/2}\Gamma \left (\frac{1-s+\delta_\eta}{2}\right )\,L(\bar \eta,1-s)
$$
where $W_\eta=i^{-\delta_\eta}N^{-1/2}\tau (\eta)$  and $\delta_\eta=\frac{1-\eta(-1)}{2}.$ 
From now until the end of this \S\,  we assume that $\eta$ is not trivial. 
For any $j\geqslant 0$ 
the special value $L(\eta,-j)$ is the algebraic integer given by
$$
L(\eta,-j)= \frac{d^j F_\eta (0)}{dt^j} \tag{25}
$$
where 
$$
F_\eta (t)=\frac{1}{\tau (\eta^{-1})}\,\underset{a\mod N}\to\sum\frac{\eta^{-1}(a)}{1-\zeta_N^ae^t}
$$
(see for example \cite{PR3, proof of Proposition 3.1.4} ). In particular
$$
L(\eta,0)=\frac{1}{\tau (\eta^{-1})}\,\underset{a\mod N}\to\sum\frac{\eta^{-1}(a)}{1-\zeta_N^a}. \tag{26}
$$
Moreover  $L(\eta,-j)=0$  if and only if $j\equiv \delta_\eta\pmod {2}.$
\newline
\newline
Let $p$ be a prime number such that $(p,N)=1.$ We fix a finite extension $L$ of $\Bbb Q_p$ containing
the values of all Dirichlet characters  $\eta$ of conductor $N$.
 The power series
$$
\scr A_{\mu_{\eta}}(X)=\,-\frac{1}{\tau (\eta^{-1})}\underset{a\mod N}\to \sum
\left ( \frac{\eta^{-1}(a)}{(1+X)\,\zeta_N^a-1}\,-\,\frac{\eta^{-1}(a)}{(1+X)^p\zeta_N^{pa}-1}\right )
$$
lies in $O_L[[X]]^{\psi=0}$ and therefore can be viewed as   the Amice transform of a unique mesure 
$\mu_\eta$ on $\Bbb Z_p^*.$  The $p$-adic $L$-functions associated to $\eta$ are defined to be 
$$
L_p(\eta\, \omega^m,s)= \int_{\Bbb Z_p^*}\omega^{m-1}(x) \bigl <x\bigr >^{-s}\mu_{\eta}(x),
\qquad 0\leqslant m\leqslant p-2.
$$
From (6) and (25) it follows  that these functions satisfy the following interpolation property
(Iwasawa theorem)
$$
 L_p(\eta\, \omega^m,1-j)= (1-(\eta \,\omega^{m-j})(p)p^{1-j})\,L(\eta\,\omega^{m-j},1-j)
 \qquad j\geqslant 1.
$$
Note that the Euler factor $1-(\eta\,\omega^{m-n})(p)p^{1-j}$ vanishes if  $m=j=1$ and $\eta (p)=1$
and that $L(\eta,0)$ does not vanish if and only if $\eta$ is odd i.e.  $\eta (-1)=-1.$
\newline
\newline
{\bf 3.2. $p$-adic representations associated to Dirichlet characters.} We continue to assume 
that $(p,N)=1.$ Set $F=\Bbb Q(\zeta_N)$,
 $G=\text{\rm Gal}(F/\Bbb Q)$ and let  $\rho\,:\,G\simeq (\Bbb Z/N\Bbb Z)^*$ denote the canonical
isomorphism normalized by $g(\zeta_N)=\zeta_N^{\rho (g)^{-1}}.$ Fix a finite extension $L/\Bbb Q_p$ containing
the values of all Dirichlet characters modulo $N$. If $\eta$ is such a character, we identify $\eta$ 
with the character $\psi\circ \rho$ of $G$ and denote by $L(\psi)$ the associated 
one-dimensional Galois representation. Let $S$ denote the set of primes
dividing $N$. 

Assume that $\eta$ is a non trivial character of conductor $N$. We need the following well known results about the Galois cohomology of $L(\eta).$

i) $H^*(\Bbb Q_l,L(\eta))=H^*(\Bbb Q_l, L(\chi\eta^{-1}))=0$ for $l\in S$.

ii)  $H^1_f(\Bbb Q,L(\eta))=0$ and  $H^1_f(\Bbb Q,L(\chi\eta^{-1}))\simeq (O_F^*\otimes_{\Bbb Z}L)^{(\eta)} .$
In particular, $H^1_f(\Bbb Q,L(\chi\eta^{-1}))=0$ if $\eta$ is odd. 

iii) The restriction of $L(\eta)$ on the decomposition group at $p$ is crystalline.
More precisely, $\Ph$ acts on $\Dc (L(\eta))$ as multiplication by $\eta (p)$ and 
the unique Hodge-Tate weight of $L(\eta)$ is $0$.

Note  that $H^0(\Bbb Q_l,L(\eta))=0$ if $l\vert N$ because in this case the inertia group acts non-trivially 
on $L(\eta)$. Together with Poincar\'e duality and the Euler characteristic formula this gives i). To prove ii)
it is enough to remark that  $H^1_f(F,\Bbb Q_p(1))\simeq   O_F^*\hat\otimes \Bbb Q_p$  (see for example \cite{Ka1,
\S5}). Finally iii) follows immediately from the definition of $\Dc$.  
\newline
\newline
Assume now that $\eta$ is odd and  $\eta (p)=1.$  Then $\Ph$ acts on $\Dc (L(\chi\eta^{-1}))$ as multiplication
by $p^{-1}$ and $D=\Dc (L(\chi\eta^{-1}))$ satisfies the conditions {\bf 1-4)} from Section 2.1.1. The isomorphism (9) writes
$$
H^1_S(\Bbb Q, L(\chi\eta^{-1}))\simeq \frac{H^1(\Bbb Q_p,L(\chi))}{H^1_f(\Bbb Q_p, L(\chi))}.
$$
{\,}
\newline
\newline
{\bf 3.3. Trivial zeros.}
\newline
{\bf 3.3.1. Cyclotomic units.} Set $F_n=F(\zeta_{p^n}).$  The collection $\bold z_{\text{\rm cycl}}=(1-\zeta_N^{p^{-n}}\zeta_{p^n})_{n\geqslant 1}$ form
a norm compartible system of units which can be viewed as an element of  $H^1_{\Iw,S}(F,L(\chi))$
using  Kummer maps $F_n^*@>>>H^1_S(F_n,L(\chi)).$ Twisting by $\ep^{-1}$ we obtain an element 
$\bold z_{\text{\rm cycl}}(-1)\in H^1_{\Iw,S}(F,L).$ Shapiro's lemma gives an isomorphism of $G$-modules 
$H^1_{\Iw}(F\otimes \Bbb Q_p,L)\simeq H^1_{\Iw}(\Bbb Q_p, L[G]^{\iota}).$ Let $\dsize e_\eta=
\frac{1}{\vert G\vert}\underset{g\in G}\to\sum \eta^{-1}(g)\,g.$ Since 
$e_\eta L[G]^{\iota}=Le_{\eta^{-1}}$ is isomorphic to $L(\eta^{-1})$
we have 
$$
e_\eta H^1_{\Iw} (F\otimes \Bbb Q_p,L)\simeq H^1_{\Iw}(\Bbb Q_p, L(\eta^{-1})).
$$
Moreover $\Dc (L[G])\simeq (L[G]\otimes F)^G\simeq L\otimes F.$ The isomorphism $\Bbb Q[G]\simeq F$
defined by $\lambda \mapsto \lambda (\zeta_N)$ induces an isomorphism $\Dc (L[G])\simeq L[G]$ and therefore
we can consider $e_\eta$ as a basis of $\Dc (L(\eta^{-1})).$
Let $\bold z_{\text{\rm cycl}}^{\eta}(-1)$ denote the image of $\bold z_{\text{\rm cycl}}(-1)$ in 
$H^1_{\Iw}(\Bbb Q_p, L(\eta^{-1})).$  We need the following properties of these elements:
\newline
\,

{\it 1) Relation to the  complex $L$-function.} Let $\bold z^{\eta^{-1}}_{\text{\rm cycl}}(-1)_0$
denote the projection of $\bold z^{\eta^{-1}}_{\text{\rm cycl}}(-1)$ on $H^1(\Bbb Q_p,L(\eta)).$ Then
$$
\exp^*_{L(\eta)}(\bold z^{\eta^{-1}}_{\text{\rm cycl}}(-1)_0)= -\left (1-\frac{\eta^{-1}(p)}{p}\right ) \,L(\eta,0)\,
e_{\eta^{-1}}
$$

{\it 2) Relation to the $p$-adic $L$-function.} Let $e^*_{\eta^{-1}}\in \Dc (L(\chi\eta^{-1}))$
be the basis which is dual to $e_{\eta^{-1}}$ and let 
$\frak L^{(\ep)}_{L(\eta),0}\,:\,H^1_{\Iw}(\Bbb Q_p,L(\eta))@>>>\scr H(\Gamma)$ denote 
the associated logarithmic map. Then
$$
\frak L^{\ep}_{L(\eta),0}(\bold z^{\eta^{-1}}_{\text{\rm cycl}}(-1))=-\bold M (\mu_\eta).
$$
We remark that 1) follows from the explicit reciprocity law of Iwasawa \cite{Iw} together with (26).
See also  \cite{Ka1, Theorem 5.12} and \cite{HK, Corollary 3.2.7}  where a more general statement is proved using the explicit reciprocity law  for 
$\Bbb Q_p(r)$. The statement 2) is a reformulation of Coleman's construction of $p$-adic $L$-functions
in terms of the large logarithmic map \cite{PR3, Proposition 3.1.4}.

\proclaim{Theorem 3.3.2}  Let $\eta$ be an odd character of conductor $N$. Assume that $p$ is a prime odd number 
such that $p\nmid N$ and $\eta (p)=1.$ Then
$$
L'(\eta\,\omega,0)\,=\, -\scr L(\eta)\, L(\eta,0)
$$ 
where $\scr L(\eta)$ is the invariant defined by (3). 
\endproclaim
\demo{Proof} It is easy to see that $\scr L(\eta)$ coincides with 
$\ell(L(\chi\eta^{-1}),D).$
Applying Proposition 2.3.4 to $V=L(\chi\eta^{-1}),$  $D=\Dc (L(\chi\eta^{-1}))$ and 
$\bold z=\bold z^{\eta^{-1}}_{\text{\rm cycl}}(-1)$    and taking into account 1-2) above
we obtain that
$$
L'_p(\eta\,\omega,0)=L'_p(\mu_z,0)=\ell(L(\chi\eta^{-1}),D) \left (1-\frac{1}{p}\right )^{-1}\left [e^*_{\eta^{-1}},\exp^*_{L(\eta)}(\bold z_0)\right ]\,=\,-\scr L(\eta)\,L(\eta,0)
$$
and the Theorem is proved.
\enddemo

\head
{\bf \S4. Trivial zeros of modular forms}
\endhead 
\flushpar 
{\bf 4.1. $p$-adic $L$-functions.}
\newline
{\bf 4.1.1. Construction of $p$-adic $L$-functions} (see \cite{AV}, \cite{Mn}, \cite{Vi}, \cite{MTT}). 
Let $f=\underset{n=1}\to{\overset\infty\to\sum} a_nq^n$ be a normalized newform  on $\Gamma_0(N)$ of weight $k$ and  character $\ep$. 
 The complex $L$-function $L(f,s)=\underset{n=1}\to{\overset \infty \to \sum}a_nn^{-s}$ decomposes into an Euler product
$$
L(f,s)=\underset{p}\to \prod E_p(f,p^{-s})^{-1}
$$
with $E_p(f,X)=1-a_pX+\ep (p)p^{k-1}X^2.$ 
Let  $p>2$ be a prime  such that the Euler factor $E_p(f,X)$ is not equal to $1$ and
let $\alpha\in \overline {\Bbb Q}_p$ be a root of the polynomial $X^2-a_pX+\ep (p)p^{k-1}.$ Assume that $\alpha$ is  not critical i.e. that   $v_p(\alpha)<k-1.$  Manin-Vishik
\cite{Mn}, \cite{Vi}, and independently Amice-Velu \cite{AV}
proved that there exists a unique distribution $\mu_{f,\alpha}$ on $\hat {\Bbb Z}^{(p)}$ 
of order $v_p(\alpha)$ such that for any Dirichlet caracter $\eta$ of conductor $M$ prime to $p$ and any
Dirichlet character $\xi$ of conductor $p^m$
$$
\int_{\hat{\Bbb Z}^{(p)}} \eta (x) \xi (x) x^{j-1}\mu_{f,\alpha}(x)=\cases
\dsize\left (1-\frac{\bar \eta (p) p^{j-1}}{\alpha} \right )\, \left (1-\frac{\beta \eta (p)}{p^j}\right ) \widetilde L(f,\eta,j)
&\text{if $1\leqslant j\leqslant k-1$ and $m=0$,}\\
\dsize\frac{p^{mj}\,\bar \eta(p^m)}{\alpha^m \tau (\bar\xi)}\,\widetilde L(f,\eta \xi^{-1},j)
&\text{if $1\leqslant j\leqslant k-1$ and $m\geqslant 1$}
\endcases
$$
where $\tau (\bar\xi)=\underset{a=1}\to{\overset{p^m-1}\to\sum} \bar\xi (a)\zeta_{p^m}^a$ 
and $\widetilde L(f,\eta,j)$ is the algebraic part of  $L(f,\eta,j)$ (see (1)).
 For us it will be more convenient to work with the distribution 
$\lambda_{f,\alpha}=x^{-1}\mu_{f,\alpha}.$ The  $p$-adic $L$-functions associated to 
$\eta \,:\,(\Bbb Z/M\Bbb Z)^*@>>> \overline{\Bbb Q}^{\,*}_p$ are defined by 
\footnote{ Our $L_{p,\alpha}(f,\eta \omega^{m},s)$  coincides with
$L_p(f,\alpha,\bar \eta \omega^{m-1},s-1)$ of \cite{MTT}}
$$
L_{p,\alpha}(f,\eta\omega^m,s)=\int_{\hat{\Bbb Z}^{(p)}} \eta   \omega^{m}(x) \bigl <x \bigr >^s \lambda_{f,\alpha} (x)
\qquad 0\leqslant m\leqslant p-2. \tag{27}
$$ 
It is easy to see that $L_{p,\alpha}(f,\eta \omega^m,s)$ is
a $p$-adic analytic function  which satisfies the following interpolation
property 
$$
L_{p,\alpha}(f,\eta \omega^m,j)= \Cal E_\alpha (f,\eta \omega^m,j)\,\,\widetilde L(f,\eta \omega^{j-m},j), \qquad 1\leqslant j
\leqslant k-1 \tag{28}
$$
where
$$
\Cal E_\alpha (f,\eta \omega^m,j)=\cases
\dsize \left (1-\frac{\bar\eta (p)p^{j-1}}{\alpha} \right )\,\left ( 1-\frac{\eta (p)\ep (p)p^{k-j-1}}{\alpha} \right )
&\text{ if $j\equiv m\pmod{p-1}$}\\
\dsize  \frac{\bar\eta (p)p^j}{\alpha \, \tau(\omega^{j-m})}
&\text{ if $j\not\equiv m\pmod{p-1}$}. 
\endcases \tag{29}
$$
\,
\flushpar
{\bf 4.1.2. $p$-adic representations associated to modular forms.}
For each prime $p$ Deligne \cite{D1} constructed a
$p$-adic representation 
$$
\rho_{f}\,:\,\text{\rm Gal} (\overline{\Bbb Q}/\Bbb Q) @>>> \text{\rm GL}(W_{f})
$$
with coefficients in a finite extension $L$ of $\Bbb Q_p.$ 
This representation has the following properties:

i) $\det \rho_{f}$ is isomorphic to $\ep \chi^{k-1}$ where $\chi$ is the cyclotomic character.

ii) If $l\nmid Np$ then the restriction of $\rho_{f}$ on the decomposition group at $l$ is unramified
and
$$
\det (1-\Fr_lX \mid W_{f})= 1-a_lX+\ep(l)\,l^{k-1}X^2
$$
(Deligne-Langlands-Carayol theorem \cite{Ca}, \cite{La}).

iii) The restriction of $\rho_{f}$ on the decomposition group at $p$ is potentially
semistable   with Hodge-Tate weights  $(0,k-1)$ \cite{Fa1}. It is crystalline if $p\nmid N$ and semistable
non-crystalline if $p \parallel N$ and $(p,  \text{\rm cond}(\ep))=1$. 
If $p\mid N$ and  $\text{\rm ord}_p(N)=\text{\rm ord}_p(\text{\rm cond}(\ep))$ 
the restriction of  $\rho_{f}$ on the decomposition group at $p$  is potentially crystalline and 
$
\Dc (W_{f})=\bD_{\text{pcris}}(W_{f})^{\text{\rm Gal}(\overline{\Bbb Q}_p /\Bbb Q_p)}
$
is  one-dimensional. In all cases
$$
\det (1-\Ph X \mid \Dc(W_{f}))= 1-a_pX+\ep(p)\,p^{k-1}X^2
$$
(Saito theorem \cite{Sa}, see also \cite{Fa2}, \cite{Ts}).
\newline
\newline
{\bf 4.1.3. Trivial zeros} (see \cite{MTT}). We say that $L_{p,\alpha}(f,\eta\omega^m,s)$ has a trivial zero
at $s=j$ if 
$$  L(f,\eta \omega^{j-m},j)\ne 0 \quad  \text{\rm and}\quad  \Cal E_{\alpha}(f,\eta \omega^m,j)=0. 
$$
From (29) it is not difficult to deduce that this occurs in the following three cases 
\cite{MTT, \S 15}:

$\bullet$ {\it The semistable case: $p\parallel N$, $k$ is even and $(p,\text{\rm cond} (\ep))=1.$} Thus $\ep (p)=0$,
$E_p(f,X)=1-a_pX$ and $a_p$ is the unique non-critical root of $X^2-a_pX.$ The restriction of $W_f$
on the decomposition group at $p$ is semistable and the eigenvalues of $\Ph$ acting on $\Dst (W_f)$
are $\alpha =a_p$ and $\beta=p\alpha.$ The module $\Dst (W_f)$ 
has a basis $\{e_\alpha,e_\beta\}$
such that $\Ph (e_\alpha)=a_pe_\alpha$, $\Ph(e_\beta)=\beta e_\beta$ and $N(e_\beta)=e_\alpha.$
Moreover $\Dc (W_f)= Le_\alpha.$
Let $\tilde \ep$ be the primitive character associated to $\ep.$ Then 
$\tilde \ep(p)\ne 0$ and  $a_p^2=\tilde \ep (p)p^{k-2}$  \cite{Li, Theorem 3}. 
Write   $a_p=\xi p^{k/2-1}$ where $\xi$ is a root of unity.  Then  $\Cal E_\alpha (f,\eta \omega^m,j)=0$  
if and only if   $j=k/2,$ $m\equiv k/2\mod{(p-1)}$ and $\bar \eta (p)= \xi.$
Therefore the $p$-adic $L$-function $L_{p,\alpha}(f,\eta \omega^{k/2},s)$ has a trivial zero at the central 
point $s=k/2$ if and only if $\bar \eta (p)=\xi.$  

$\bullet${\it The crystalline case: $p\nmid N$.} The restriction of $W_f$ on the decomposition group at $p$ is crystalline 
and by Deligne \cite{D2} one has  $\vert \alpha\vert=p^{(k-1)/2}.$ 
Write $\alpha= \xi p^{\frac{k-1}{2}}$ with $\vert \xi\vert =1.$ 
Then $\Cal E_\alpha (f,\eta \omega^m,j)$ vanishes if and only if $m\equiv j\pmod{p-1}$,  $k$ is odd, and either
$j=\dsize \frac{k+1}{2}$ and $\bar \eta (p)=\xi$ or $j=\dsize \frac{k-1}{2}$ and $\eta (p)\ep (p)=\xi.$
The $p$-adic $L$-function $L_{p,\alpha}(f,\eta \omega^{\frac{k+1}{2}},s)$ has a trivial zero at the near-central point
$s=\frac{k+1}{2}$ if and only if $\alpha=\bar \eta (p)\,p^{\frac{k-1}{2}}$ and $L_{p,\alpha}(f,\eta \omega^{\frac{k-1}{2}},s)$ has a trivial zero at $s=\frac{k-1}{2}$ if and only if $\alpha=\eta (p)\,\ep (p)\, p^{\frac{k-1}{2}}.$
    
$\bullet${\it The potentially crystalline case: $p\mid N$ and $\text{\rm ord}_p(N)=\text{\rm ord}_p(\text{\rm cond} (\ep)).$} 
 One has $E_p(f,X)=1-a_pX$ and $\alpha=a_p$ is the unique non-critical root
of $X^2-a_pX.$ Moreover  $\tilde \ep (p)=0$ and
it can be shown that  $|a_p|=p^{\frac{k-1}{2}}$ \cite{O}, \cite{Li}. 
The restriction of $W_f$ on the decomposition group at $p$ is potentially crystalline and $\Dc (W_f)$ is
a one-dimensional vector space on which $\Ph$ acts as multiplication by $a_p.$
The factor $\Cal E_\alpha (f,\eta \omega^m,j)$ vanishes if and only if $k$ is odd, $j=m=\frac{k+1}{2}$ and $a_p=\bar \eta (p)p^{\frac{k-1}{2}}.$
The  $p$-adic $L$-function $L_{\alpha,p}(f,\eta \omega^{\frac{k+1}{2}},s)$ has a trivial zero at the near-central point $s=\frac{k+1}{2}$ if and only if $a_p=\eta (p)p^{\frac{k-1}{2}}.$
\newline
\newline
If $\eta$ is a Dirichlet character of conductor $M$,  
the twisted modular form $f_{\eta}=\underset{n=1}\to{\overset\infty\to\sum}\eta (n)\,a_nq^n$ is not 
necessarily primitive, but there exists a unique normalized newform $f\otimes\eta$ such that
$$
L(f,\eta,s)=L(f\otimes \eta,s)\,\,\underset{l\mid M}\to \prod  E_l (f\otimes\eta,l^{-s}) 
$$ 
(see for example \cite{AL}). 
Write $L(f\otimes \eta,s)=\underset{n=1}\to{\overset \infty \to \sum}\dsize\frac{a_{\eta,n}}{n^s}.$
If $p\nmid M$, the Euler factors at $p$ of $L_M(f\otimes\eta,s)$  and $L(f,\eta,s)$  coincide and $\alpha_\eta=\alpha \eta (p)$
is a root of $X^2-a_{\eta,p}X+\ep (p)\eta^2(p)p^{k-1}.$ It is easy to see that $\Cal E_{\alpha_\eta}(f_\eta, \omega^{m},j)=
\Cal E_\alpha (f,\eta\omega^m,j)$ and 
from the interpolation formula (28)  it follows immediately that
the behavior of $L_{p,\alpha}(f,\eta\omega^{m},s)$ and $L_{p,\alpha_\eta}(f\otimes \eta,\omega^m,s)$ is essentially the same. 
Therefore the general case reduces 
to the case of the trivial character $\eta.$ 
\newline
\newline
{\bf 4.2. Selmer groups and  $\ell$-invariants of modular forms.} 
\newline
{\bf 4.2.1. The Selmer group.} From now until the end of this $\S$  we assume that $L_{p,\alpha}(f,\omega^{k_0},s)$ has a trivial zero at $k_0.$ Thus $k_0=k/2$ in the semistable case and $k_0=\frac{k\pm 1}{2}$ in the crystalline or potentially crystalline case. Set $V_f=W_f\left (k_0\right ).$
Let $f^*$ denote the complex conjugation of $f$ i.e. $f^*=\underset{n=1}\to{\overset \infty \to \sum}
\bar a_nq^n.$  The canonical pairing
$W_f\times W_{f^*}@>>>L(1-k)$ induces an isomorphism $W_{f^*}\left (k-k_0\right )\simeq V_f^*(1).$
We need the following  basic results
about the Galois cohomology of $V_f$: 

i) $H^0(\Bbb Q_p,V_f)=H^0(\Bbb Q_p,V_f^*(1))=0$ and $\dim_L H^1(\Bbb Q_p,V_f)=\dim_L H^1(\Bbb Q_p,V_f^*(1))=2.$ 

ii) $H^1_f(\Bbb Q_p,V_f)$ and $H^1_f(\Bbb Q_p,V_f^*(1))$ are  one-dimensional $L$-vector spaces.

iii) $H^1_f(\Bbb Q,V_f)=H^1_f(\Bbb Q,V_f^*(1))=0.$

We remark that using the fact that the Hodge-Tate weights of $W_f$ are $0$ and $k-1$ 
and the eigenvalues of $\Ph$ on $\Dc (W_f)$ have absolute value $p^{(k-1)/2}$ (respectively $p^{k/2}$)
in the crystalline and potentially crystalline case (respectively in the semistable case) one  deduce
that $H^0(\Bbb Q_p,W_f(m))=0$ for all $1\leqslant m\leqslant k-1$ (see \cite{Ka2, Proposition 14.12 and
Section 13.3}). Applying Poincar\'e duality and  the Euler characteristic formula we obtain i). Next ii) follows from i) together with  the formula
$$
\dim_L H^1_f(\Bbb Q_p,V_f)=\dim_L t_{V_f}(L)+\dim_L H^0(\Bbb Q_p,V_f).
$$
Finally iii) is a deep result of Kato \cite{Ka2, Theorem 14.2}. Note that in the semistable case 
we assume that $L(f,k/2)\ne 0.$

From i-iii) above it follows that $V_f$ satisfies the conditions {\bf 1-3)} of Section 2.2.1.
Assume that $k_0\geqslant \frac{k+1}{2}.$  This holds
automatically in the semistable  ($k_0=k/2$) and potentially crystalline  ($k_0=\frac{k+1}{2}$) cases. In the crystalline case
$\alpha^*=\ep^{-1}(p)\alpha$ is a root of $1-\bar a_pX +\ep^{-1}(p)p^{k-1}X^2$ and using the functional equation
for $p$-adic $L$-functions one can reduce the study of   $L_{p,\alpha}(f,\omega^{k_0},s)$ at $s=\frac{k-1}{2}$ 
to the study of $L_{p,\alpha^*}(f^*,\omega^{k_0+1},s)$  at $s=\frac{k+1}{2}.$

\proclaim{Lemma 4.2.2} Assume that $L_{p,\alpha}(f,\omega^{k_0},s)$ has a trivial zero 
on the right of the central point (i.e.  $k_0\geqslant k/2$). Then $D_\alpha=\Dc (V_f)^{\Ph=p^{-1}}$ is a
one dimentional $L$-vector space which satisfies one of the conditions a-b) of Proposition 2.1.2.
\endproclaim
\demo{Proof}
From 4.1.2 it follows that $\alpha p^{-k_0}=p^{-1}$ is an eigenvalue of $\Ph$ acting on $\Dc (V_f).$
Thus $\dim_L D_\alpha \geqslant 1.$ If $\dim_L D_\alpha=2$ then $V_f$ would be crystalline 
and $\Ph$ would act on $\Dc (V_f)$ as multiplication by $p^{-1}.$ This contadicts the weak
admissibility of $\Dc (V_f).$  Finally $D_\alpha$ satisfies a) in the crystalline and potentially crystalline cases
and b) in the semistable case.
\enddemo

{\bf 4.2.3. The $\ell$-invariant of modular forms.} From Lemma 4.2.2 it follows that if $L_{p,\alpha}(f,\omega^{k_0},s)$ has a trivial zero
at $k_0\geqslant k/2$ the $\ell$-invariant  $\ell (V_f,D_\alpha)$ is well defined. To simplify notation we will denote it 
by $\ell_\alpha (f).$ The general definition of the $\ell$-invariant can be made more explicit in the case 
of modular forms.

{\it $\bullet$ The semistable case.} Let $\{e_\alpha,e_\beta\}$ denote the basis of $\Dst (W_f)$  as in 4.1.3. 
In \cite{Ben2, Proposition 2.3.7} it is proved that
$$
\ell_\alpha (f)=\scr L_{\text{FM}}(f) \tag{30}
$$
where $\scr L_{\text{FM}}(f)$ is the Fontaine-Mazur invariant \cite{Mr} which is defined as the unique element of $L$
such that
$$
e_\beta + \scr L_{\text{FM}}(f)\, e_\alpha \in \F^{k-1}\Dst (W_f).
$$

{\it $\bullet$ The crystalline and potentially crystalline cases.}  
The $(\Ph,\Gamma)$-module  $\Ddagrig (V_f)\cap \bigl(D_{\alpha}
\otimes_L \CR_L[1/t]\bigr )$  is isomorphic
to $\CR_L(\delta)$ with $\delta (x)=|x|x^{\frac{k+1}{2}}$  and the exact sequence (11) writes 
$$
0@>>> \CR_L(\delta)@>>>\Ddagrig (V_f) @>>> \CR_L(\delta')@>>>0 
$$ 
for some character $\delta' \,:\,\Bbb Q_p^*@>>>L^*$. 
Since 
$\dim_L H^1(\CR_L(\delta))=2$  we have 
$
H^1(\Bbb Q_p,V_f)\simeq H^1(\CR_L(\delta)). 
$   
Therefore  $H^1_{D_\alpha}(V_f)= H^1_{f,\{p\}}(\Bbb Q,V_f)$
and $\scr L_\alpha(f)=\scr L(V_f,D_\alpha)$ is the slope of the image of the localization map
$
H^1_{f,\{p\}}(\Bbb Q,V_f) @>>>H^1(\Bbb Q_p,V_f)
$
under the canonical decomposition (6) 
$$H^1(\Bbb Q_p,V_f)\simeq H^1_f(\CR_L(\delta))\times H^1_c(\CR_L(\delta)).
$$
The formula (19) writes 
$$
\ell_\alpha (f)=-\scr L_\alpha (f). \tag{31}
$$
\newline
\newline
{\bf 4.3. The main result.}
\newline
{\bf 4.3.1. Kato's Euler systems.} 
Using the theory of modular units Kato \cite{Ka2} constructed an element $\bold z_{\text{Kato}}\in H^1_{\Iw,S}(W_{f^*})$ 
which is closely related to the complex and the $p$-adic $L$-functions via the Bloch-Kato exponential map. The CM-case
was considered before by Rubin \cite{Ru}. Set 
$$\bold z_{\text{Kato}}(j)=\text{Tw}^\ep_{j}(\bold z_{\text{Kato}})\in H^1_{\Iw,S}(W_{f^*}(j))
$$
and denote by ${\bold z_{\text{Kato}}(j)}_0=\text{\rm pr}_{0}(\bold z_{\text{Kato}}(j))$ the projection of 
$\bold z_{\text{Kato}}(j)$ on $H^1_S(W_{f^*}(j))$. 
The following statements are direct analogues of properties 1-2) of cyclotomic units from Section 3.3.1:
\newline
\,

{\it 1) Relation to the complex $L$-function.} One has
$$
\exp^*_{W_{f^*}(j)} ( {\bold z_{\text{Kato}} (j)}_0)= \Gamma (k-j)^{-1}E_p\left (f,p^{k-j}\right)\, \widetilde L\left(f,k-j\right ) 
 \omega^*_j, \qquad 1\leqslant j\leqslant k-1 \tag{32}
$$
for some canonical basis $\omega^*_j$ of $\F^0\Dc (W_f^*(j))$
\cite{Ka2, Theorem 12.5}. Note that $\omega^*_{j+1}=\omega^*_{j}\otimes e_1$ where $e_1=\ep^{-1}\otimes t$.

{\it 2) Relation to the $p$-adic $L$-function.} Fix a generator $d_\alpha$ of $D_\alpha.$ Let $\Log^{(\alpha),\ep}_{W_{f^*}(k),1}$ denote the large logarithmic map
$\Log^{(\ep)}_{W_{f^*}(k),1,\eta}$
associated to $\eta=d_\alpha \otimes e_{\frac{k+1}{2}}  \in \Dc (W_f).$
Then
$$
\Log^{(\alpha),\ep}_{W_{f^*}(k),1}(\bold z_{\text{\rm Kato}}(k))=
\text{\bf M}(\lambda_{f,\alpha})\, \left [d_\alpha \otimes e_{\frac{k+1}{2}}\,,\,\omega^*_k \right ]_{W_f} \tag{33}
$$
\cite{Ka2, Theorem 16.2}.

We can now prove the main result of this paper.

\proclaim{Theorem 4.3.2} Let $f$ be a newform on $\Gamma_0(N)$ of character $\ep$ and  weight $k$ and
let  $p$ be an odd prime.
Assume that  the $p$-adic $L$-function 
$L_{p,\alpha}\left (f,\omega^{k_0},s \right )$ has a trivial zero at $s=k_0\geqslant k/2$.
Then 
$$
L'_{p,\alpha}\left (f, \omega^{k_0},k_0\right )=\ell_\alpha (f)\,\left (1-\frac{\ep (p)}{p}
\right ) \widetilde L \left (f,k_0\right ).
$$
\endproclaim 
\demo{Proof} To simplify notation set   $\bold z=\bold z_{\text{\rm Kato}}\left (k-k_0\right )$. 
By Lemma 1.3.4 one has  
$$\Log_{V^*_f(1),1-k_0}^{(\alpha),\ep}(\bold z)=\Tw_{k_0}\left ( 
\Log_{W_{f^*}(k),1}^{(\alpha),\ep}(\bold z_{\text{\rm Kato}}(k))\right ).
$$
Let $\mu_{\bold z}$ be the distribution defined by $\text{\bf M}(\mu_{\bold z})= 
\Log_{V^*_f(1),1-k_0}^{(\alpha),\ep}(\bold z).$ 
Then (33) gives
$$
\text{\bf M}(\mu_{\bold z})=\Tw_{k_0}\left (\text{\bf M} (\lambda_{f,\alpha})\right )\, 
{[d_\alpha \otimes e_{k_0},\omega_k^*]}_{W_f}\,=\,\Tw_{k_0}\left (\text{\bf M} (\lambda_{f,\alpha})\right )\,  {[d_\alpha ,\omega_{k_0-1}^*]}_{V_f} 
$$
and from (8) and (27) it follows that
$$
L_p(\mu_{\bold z},s)= L_{p,\alpha}\left (f,\omega^{k_0},s+k_0 \right )\, {\left [d_\alpha ,\omega_{k_0-1}^* \right ]}_{V_f}.
$$ 
Now, applying  Proposition 2.2.4 we obtain
$$
 L'_{p,\alpha}\left (f,\omega^{k_0},k_0 \right )\, {\left [d_\alpha ,\omega_{k_0-1}^* \right ]}_{V_f}\,=\,
 \ell_{\alpha}(f)\,\Gamma \left (k_0\right )\,\left (1-\frac{1}{p}\right )^{-1}\,
 {\left [d_\alpha , \exp^*_{V^*_f(1)}(\bold z_0) \right ]}_{V_f}. \tag{34}
$$
On the other hand, for $j=k-k_0$ the formula (32) gives
$$
\exp^*_{V_f^*(1)} (\bold z_0 ) \,=\,\Gamma \left (k_0\right )^{-1}E_p \left (f,p^{k_0}\right )\,
\widetilde L \left (f,k_0\right )\,\omega^*_{k_0-1}. \tag{35}
$$
Since $E_p \left (f,p^{k_0}\right )\,=\,\dsize \left (1-\frac{1}{p} \right )\,
\left (1-\frac{\ep (p)}{p}\right )$ and $\left [d_\alpha ,\omega_{k_0-1}^* \right ]_{V_f}\ne 0,$ from (34) and (35)
we obtain that
$$
L'_{p,\alpha}\left (f, \omega^{k_0},k_0\right )=\ell_{\alpha} (f)\,\left (1-\frac{\ep (p)}{p}
\right ) \widetilde L \left (f,k_0\right )
$$
and the Theorem is proved.
\enddemo

\flushpar
\newline
{\bf Corollaries 4.3.3.} 1) In the semistable case $k$ is even and $\ep (p)=0$. Theorem 4.3.2 together with (30)
give the Mazur-Tate-Teitelbaum conjecture 
$$
L'_{p,\alpha} (f, \omega^{k/2},k/2 )=\scr L_{\text{\rm FM}} (f)\,\widetilde L \left (f,k/2\right )
$$
and our proof can be seen as a revisiting of Kato-Kurihara-Tsuji approach using the theory
of $(\Ph,\Gamma)$-modules.

2) In the crystalline and potentially crystalline cases Theorem 4.3.2 writes
$$
L'_{p,\alpha}\left (f, \omega^{\frac{k+1}{2}},\frac{k+1}{2}\right )=-\scr L_\alpha (f)\,\left (1-\frac{\ep (p)}{p}
\right ) \widetilde L \left (f,\frac{k+1}{2}\right )
$$
(see (31)).

\enddocument 

\Refs\nofrills{References} \widestnumber \key{MTT} 
\ref
\key{AL}
\by A.O.L. Atkin and W. Li
\paper Twists of Newforms and Pseudo-Eigenvalues of $W$-Operators
\jour Invent. Math.
\vol 48
\yr 1978
\pages 221-243
\endref 
\ref \key AV
\by Y. Amice and J. V\'elu
\paper Distributions $p$-adiques associ\'ees aux s\'eries de Hecke
\jour Ast\'erisque
\vol 24-25
\yr 1975
\pages 119-131
\endref
\ref{}
\key{Bel}
\by J. Bella\"{\i}che
\paper Critical $p$-adic $L$-functions
\jour Inventiones Math. 
\yr 2012
\vol 189
\pages 1-60  
\endref
\ref
\key{BC} \by
J. Bella\"{\i}che,  G. Chenevier
\paper $p$-adic families of Galois representations and higher rank Selmer groups
\jour Ast\'erisque
\vol 324
\yr 2009
\pages 314 pages
\endref
\ref{}\key{Ben1}
\by D. Benois
\paper On Iwasawa theory of crystalline representations
\jour Duke Math. J.
\yr 2000
\vol 104
\pages 211-267
\endref
\ref{}
\key{Ben2} 
\by D. Benois
\paper A generalization of Greenberg's $\scr L$-invariant
\jour Amer. J. of Math. 
\vol 133
\issue 6
\yr 2011
\pages 1573-1632
\endref
\ref{}
\key{Ben3}
\by D. Benois
\paper Infinitesimal deformations and the $\ell$-invariant
\jour Documenta Math. Extra Volume: Andrei A. Suslin's Sixtieth Birthday
\yr 2010
\pages 5-31 
\endref
\ref \key Ber1 \by L. Berger \paper Repr\'esentations $p$-adiques
et \'equations diff\'erentielles \jour Invent. Math. \yr 2002 \vol 148 \issue 2 \pages
219-284
\endref
\ref
\key Ber2
\by L. Berger
\paper Bloch and Kato's exponential map: three explicit
formulas
\jour Doc. Math., Extra Vol.
\yr 2003
\pages 99-129
\endref
\ref
\key{Ber3}
\by L. Berger
\paper Limites de repr\'esentations cristallines
\jour Compos. Math.
\yr 2004
\vol 140
\issue 6
\pages 1473-1498
\endref
\ref \key Ber4 \by L. Berger \paper Equations diff\'erentielles $p$-adiques et $(\Ph,N)$-modules filtr\'es
\jour Ast\'erisque 
\vol 319
\yr 2008
\pages 13-38 
\endref
\ref \key {BK} \by S. Bloch, K. Kato \paper $L$-functions and
Tamagawa numbers of motives \jour Grothendieck Fest-\linebreak
schrift, vol. 1 \pages 333-400 \yr 1990
\endref
\ref
\key{Ca}
\by H. Carayol
\paper Sur les repr\'esentations $l$-adiques associ\'ees aux formes modulaires de Hilbert
\jour Ann. Sci. ENS
\yr 1986
\vol 19
\pages 409-468
\endref
\ref
\key{CC1}
\by F. Cherbonnier and P. Colmez
\paper Repr\'esentations $p$-adiques surconvergentes
\jour Invent. Math.
\yr 1998
\vol 133
\pages 581-611
\endref
\ref
\key{CC2}
\by F. Cherbonnier and P. Colmez
\paper Th\'eorie d'Iwasawa des repr\'esentations $p$-adiques d'un corps local
\jour J. Amer. Math. Soc.
\vol 12
\yr 1999
\pages 241-268
\endref
\ref{}
\key{Co}
\by R. Coleman
\paper A $p$-adic Shimura isomorphism and $p$-adic periods of modular forms
\jour Contemp. Math.
\vol 165
\yr 1994
\pages 21-51
\endref
\ref{}
\key{CI}
\by R. Coleman and A. Iovita
\paper Hidden structures on semi-stable curves
\jour Ast\'erisque
\vol 331
\yr 2010
\pages 179-254
\endref
\ref{}
\key{Cz1}
\by P. Colmez
\paper  Th\'eorie d'Iwasawa des repr\'esentations de de Rham d'un corps local
\jour Ann. of Math.
\yr 1998
\vol 148
\pages 485-571
\endref
\ref
\key{Cz2}
\by P. Colmez 
\paper Repr\'esentations crystallines et repr\'esentations de hauteur finie
\jour J. Reine und Angew. Math.
\yr 1999
\vol 515
\pages 119-143
\endref
\ref
\key{Cz3}
\by P. Colmez
\paper La conjecture de Birch et Swinnerton-Dyer $p$-adique
\jour S\'eminaire Bourbaki 2002/03, Ast\'erisque
\vol 294
\yr 2004
\pages 251-319
\endref
\ref{}
\key{Cz4}
\by P. Colmez 
\paper Z\'eros suppl\'ementaires de fonctions $L$ $p$-adiques
de formes modulaires
\jour Algebra and Number Theory,
Hindustan book agency 2005
\pages 193-210
\endref
\ref
\key{Cz5}
\by P. Colmez
\paper Repr\'esentations triangulines de dimension $2$
\jour Ast\'erisque
\vol 319
\yr 2008
\pages 213-258
\endref
\ref{}
\key{Cz6}
\by P. Colmez 
\paper Fonctions d'une variable $p$-adique
\jour Ast\'erisque
\vol 330
\yr 2010
\pages 13-59
\endref 
\ref{}
\key{DDP}
\by S. Dasgupta, H. Darmon and R. Pollack
\paper Hilbert modular forms and the Gross-Stark conjecture
\jour Annals of Math. (to appear)
\vol
\yr
\pages
\endref
\ref{}\key{D1}
\by P. Deligne
\paper Formes modulaires et repr\'esentations $l$-adiques
\jour S\'em. Bourbaki 1968/69, Lecture Notes in Math.
\vol 179
\yr 1971
\pages 139-172
\endref
\ref 
\key{D2}
\by P. Deligne
\paper La conjecture de Weil I
\jour Publ. Math. de l'IHES
\vol 43
\yr 1974
\pages 273-307
\endref
\ref\key{FG}
\by B. Ferrero and R. Greenberg
\paper On the behavior of $p$-adic $L$-functions at $s=0$
\jour Invent. Math.
\vol 50
\yr 1978/79
\pages 91-102
\endref
\ref{}
\key{Fa1}
\by G. Faltings
\paper Hodge-Tate structures and modular forms
\jour Math. Ann.
\yr 1987
\vol 278
\pages 133-149
\endref
\ref{}
\key{Fa2}
\by G. Faltings 
\paper Crystalline cohomology and $p$-adic Galois representations
\jour in "Algebraic Analysis, Geometry and Number Theory", 
John Hopkins University Press, Baltimore MD, 1989
\pages 25-80
\endref
\ref \key{Fo1} \by J.-M. Fontaine \paper Repr\'esentations
$p$-adiques des corps locaux
\jour  in "The Grothendieck
Festschrift", vol. II, Birkh\"auser, Boston \yr 1991 \pages
249-309
\endref
\ref \key{Fo2} \by J.-M. Fontaine \paper Le corps des p\'eriodes
$p$-adiques \jour Ast\'erisque \vol 223 \yr 1994 \pages 59-102
\endref
\ref \key{Fo3} \by J.-M. Fontaine \paper Repr\'esentations
$p$-adiques semi-stables \jour Ast\'erisque \vol 223 \yr 1994
\pages 113-184
\endref
\ref
\key{FP}
\by J.-M. Fontaine and B. Perrin-Riou
\paper Autour des conjectures de Bloch et Kato; cohomologie
galoisienne et valeurs de fonctions $L$
\jour in "Motives", Proc. Symp. in Pure Math.,
vol. 55, part 1,
\yr 1994
\pages 599-706
\endref
\ref\key{Gre1}
\by R. Greenberg
\paper Iwasawa theory of $p$-adic representations
\jour Adv. Stud. in Pure Math. 
\vol 17
\yr 1989
\pages 97-137
\endref
\ref \key{Gre2}
\by R. Greenberg
\paper Trivial zeros of $p$-adic $L$-functions
\jour Contemp. Math.
\vol 165
\yr 1994
\pages 149-174
\endref
\ref \key{GS}
\by R. Greenberg and G. Stevens
\paper $p$-adic $L$-functions and $p$-adic periods of modular forms
\jour Invent. Math.
\yr 1993
\vol 111
\pages 407-447
\endref
\ref
\key{Gro}
\by B. Gross
\paper $p$-adic $L$-series at $s=0$
\jour J. Fac. Sci. Univ. Tokyo, sect. IA Math.
\yr 1981
\vol 28
\pages 979-994
\endref
\ref{}\key{GK}
\by B. Gross and N. Koblitz
\paper Gauss sums and the $p$-adic $\Gamma$-function
\jour Annals of Math.
\vol 109
\yr 1979
\pages 569-581
\endref
\ref \key H1 \by L. Herr \paper  Sur la cohomologie galoisienne
des corps $p$-adiques \jour Bull. Soc. Math. France \yr 1998 \vol
126 \pages 563-600
\endref
\ref \key {H2}
\by L. Herr
\paper Une approche nouvelle de la dualit\'e locale de Tate
\jour Math. Annalen
\yr 2001
\vol 320
\pages 307-337
\endref
\ref{}
\key{HK}
\by A. Huber and G. Kings
\paper Bloch-Kato conjecture and  main conjecture of Iwasawa theory for Dirichlet characters
\jour Duke Math. J.
\vol 119
\yr 2003
\pages 393-464
\endref
\ref{}
\key{Iw}
\by K. Iwasawa
\paper On explicit formulas for the norm residue symbol
\jour J. Math. Soc. Japan
\vol 20
\yr 1968
\pages 151-165 
\endref 
\ref \key{JS}
\by H. Jacquet and J. Shalika
\paper A non-vanishing theorem for zeta functions on $\text{\rm GL}_n$
\jour Invent. Math.
\vol 38
\yr 1976
\pages 1-16
\endref
\ref{}\key{Ka1}
\by K. Kato
\paper Iwasawa theory and $p$-adic Hodge theory
\jour Kodai Math. J.
\vol 16
\yr 1993
\pages 1-31
\endref
\ref{}
\key{Ka2}
\by K. Kato
\paper $p$-adic Hodge theory and values of zeta-functions of modular forms
\jour Asterisque
\vol 295
\yr 2004
\pages 117-290
\endref
\ref{}
\key{KKT}
\by K. Kato, M. Kurihara and T. Tsuji
\paper Local Iwasawa theory of Perrin-Riou and syntomic complexes, preprint 
\yr 1996
\endref  
\ref \key Ke \by K. Kedlaya \paper A p-adic monodromy theorem
\jour Ann. of Math. \vol 160 \yr 2004 \pages 93-184
\endref
\ref{}
\key{La}
\by R.P. Langlands
\paper Modular forms and $l$-adic representations
\jour in "Modular forms of one variable II", Lecture Notes in Math.
\vol 349 
\pages 361-500
\endref
\ref{}\key {Li}
\by W. Li
\paper Newforms and Functional Equations
\jour Math. Ann.
\yr 1975
\vol 212
\pages 285-315
\endref
\ref{}
\key{Liu}
\by R. Liu
\paper Cohomology and duality for $(\Ph,\Gamma)$-modules over the Robba ring
\jour Int. Math. Research Notices
\vol 3
\yr 2008
\pages Art. ID rnm150,  32 pages
\endref
\ref{}
\key{Mn}
\by Y. Manin
\paper Periods of cusp forms and $p$-adic Hecke series
\jour Math. USSR Sbornik 
\vol 92
\yr 1973
\pages 371-393
\endref
\ref
\key{Mr}
\by B. Mazur
\paper On monodromy invariants occuring in global arithmetic and Fontaine's theory
\jour Contemp. Math.
\vol 165
\yr 1994
\pages 1-20
\endref
\ref{}\key{MTT}
\by B. Mazur, J.Tate and J. Teitelbaum
\paper On $p$-adic analogues of the conjectures of Birch and Swinnerton-Dyer
\jour Invent. Math.
\yr 1986
\vol 84
\pages 1-48
\endref
\ref 
\key{Mi} 
\by T. Miyake
\paper Modular forms
\jour Springer Monographs in Mathematics, Springer, 2006
\pages 335 pages
\endref
\ref{}\key{Ne}
\by J. Nekov\'a\v{r}
\paper On $p$-adic height pairing
\jour S\'eminaire de Th\'eorie des Nombres, Paris 1990/91,
Progress in Math.
\vol 108
\yr 1993
\pages 127-202
\endref
\ref{}\key{O}
\by A. Ogg
\paper On the eigenvalues of Hecke operators
\jour Math. Ann.
\yr 1969
\vol 179
\pages 101-108
\endref
\ref{}\key{Or}
\by L. Orton
\paper On exceptional zeros of $p$-adic $L$-functions associated
to modular forms
\jour PhD thesis, University of Nottingham
\yr 2003
\endref
\ref
\key{PR1}
\by B. Perrin-Riou
\paper Th\'eorie d'Iwasawa et hauteurs $p$-adiques
\jour Invent. Math.\vol 109\yr 1992\pages 137-185\endref
\ref 
\key PR2
\by B. Perrin-Riou
\paper Th\'eorie d'Iwasawa des repr\'esentations $p$-adiques sur un corps local
\jour Invent. Math.\vol 115\yr 1994\pages 81-149
\endref
\ref{}
\key{PR3}
\by B. Perrin-Riou
\paper La fonction $L$ $p$-adique de Kubota-Leopoldt
\jour Contemp. Math.
\vol 174
\yr 1994
\pages 65-93
\endref
\ref{}\key{PR4}
\by B. Perrin-Riou
\paper Z\'eros triviaux des fonctions $L$ $p$-adiques
\jour Compositio Math.
\yr 1998
\vol 114
\pages 37-76
\endref
\ref{}
\key{PR5}
\by B. Perrin-Riou
\paper Quelques remarques sur la th\'eorie d'Iwasawa des courbes elliptiques
\jour  Number Theory for Millennium, III (Urbana, IL, 2000)
\pages 119-147
\endref
\ref{}\key{PS}
\by R. Pollack and G. Stevens
\paper Critical slope  $p$-adic $L$-functions
\jour  To appear in the "Journal of the London Math. Soc"
\pages Preprint  available on 
 http://math.bu.edu/people/rpollack
\endref
\ref
\key{Pt}
\by J. Pottharst
\paper Cyclotomic Iwasawa theory of motives
\jour Preprint
\yr 2012
\pages 25 pages
\endref
\ref{}\key{Ru}
\by K. Rubin
\paper Elliptic curves with complex multiplication and the conjecture of Birch
and Swinnerton-Dyer
\jour Invent. Math.
\vol 64
\yr 1981
\pages 455-470
\endref
\ref{}
\key{Sa}
\by T. Saito
\paper Modular forms and $p$-adic Hodge theory
\jour Invent. Math.
\yr 1997
\vol 129
\pages 607-620
\endref
\ref{}\key{St}
\by G. Stevens
\paper Coleman's $\scr L$-invariant and families of modular forms 
\jour Ast\'erisque 
\yr 2010
\vol 311
\pages 1-12
\endref
\ref{}
\key{Ts}
\by T. Tsuji
\paper $p$-adic etale cohomology and crystalline cohomology in the semistable reduction case
\jour Invent. Math.
\yr 1997
\vol 137
\pages 233-411
\endref
\ref{}
\key{Vi}
\by M. Vishik
\paper Non archimedean measures connected with Dirichlet series
\jour Math. USSR Sbornik
\vol 28
\yr 1976
\pages 216-228
\endref
\ref
\key{Wa}
\by N. Wach
\paper Repr\'esentations $p$-adiques potentiellement cristallines
\jour Bull. Soc. Math. France 
\vol 124
\yr 1996 
\pages 375-400
\endref
\endRefs
\bye